\documentclass[11pt]{amsart}
\usepackage{latexsym,amssymb,amsmath,amscd}
\numberwithin{equation}{section}
\theoremstyle{plain}

\newtheorem{definition}{{\bf DEFINITION}}[section]
\newtheorem{corollary}{{\bf COROLLARY}}[section]
\newtheorem{example}{{\bf EXAMPLE}}[section]
\newtheorem{proposition}{{\bf PROPOSITION}}[section]
\newcommand{\bq}{\begin{equation}}
\newcommand{\bea}{\begin{array}}
\newcommand{\eea}{\end{array}}

\newcommand{\ga}{\alpha}
\newcommand{\gep}{\epsilon}
\newcommand{\gD}{\Delta}
\newcommand{\gl}{\lambda}
\newcommand{\gL}{\Lambda}
\newcommand{\gb}{\beta}
\newcommand{\ot}{\otimes}
\newcommand{\mf}{\mathfrak}
\newcommand{\mc}{\mathcal}

\newcommand{\wg}{\wedge}

\newcommand{\lcrb}{>\!\!\vartriangleleft\!\!\!\!\cdot}

\newcommand{\ci}{\circ}
\newcommand{\ul}[1]{\underline{#1}}
\newcommand{\ol}[1]{\overline{#1}}

\newcommand{\go}{\omega}
\newcommand{\gO}{\Omega}
\newcommand{\gG}{\Gamma}
\newcommand{\gt}{\theta}
\newcommand{\gs}{\sigma}
\newcommand{\gz}{\zeta}
\newcommand{\gag}{\gamma}
\newcommand{\gd}{\delta}
\newcommand{\pp}{\partial}

\newcommand{\tl}{\tilde}
\newcommand{\na}{\nabla}
\newcommand{\gk}{\kappa}
\newcommand{\nm}{\left[\begin{array}{c}
n\\
m\end{array}\right]}

\newcommand{\bl}{\blacklozenge}
\newcommand{\bs}{\blacksquare}

\newcommand{\gT}{\Theta}
\newcommand{\ddg}{\ddagger}

\title{REMARKS ON QUANTUM TRANSMUTATION}
\author{Robert Carroll\\University of Illinois, Urbana, IL 61801}

\date{January, 2001 - email: rcarroll@math.uiuc.edu}

\begin{document}

\bibliographystyle{plain}

\begin{abstract} {It is shown how formulas of the author for general operator
transmutation can be adapted to a quantum group context}
\end{abstract}
\maketitle




\section{INTRODUCTION}
\renewcommand{\theequation}{1.\arabic{equation}}
\setcounter{equation}{0}

We want to sketch here a program and a framework along with a few ideas about
implementing it.  First we indicate briefly some results of Koelink and Rosengren
\cite{kf}, involving tramsmutation kernels for little q-Jacobi functions.  This is
modeled in part on earlier work of Koornwinder \cite{kzz} and others on classical
Fourier-Jacobi, Abel, and Weyl transforms, plus formulas in q-hypergeometric
functions (see e.g. \cite{kd,kj,ki}).  Here transmutation refers to intertwining of
operators (i.e. $QB=BP$) and the development
gives a q-analysis version of some transforms and
intertwinings whose classical versions were embedded as examples in a general theory
of transmutation of operators
by the author (and collaborators) in \cite{bv,cw,cv,cu,cx} (see also the references
in these books and papers).  It is clear that many of the formulas from the extensive
general theory will have a q-analysis version (as in \cite{kf})
and the more interesting problem
here would seem to be that
of phrasing matters entirely in the language of quantum groups and ultimately connecting
the theory to tau functions in the spirit indicated below.  In this direction we sketch
some of the ``canonical" development of the author in \cite{cw}.  
\\[3mm]\indent
We mention also a quantum group notion of transmutation
\bq\label{500}
\boxed{Braided\,\,
stuff}\begin{array}{c}
\stackrel{Transmutation}{\longleftarrow}\\
\stackrel{Bosonization}{\longrightarrow}
\end{array}
\boxed{Quantum\,\, stuff}
\end{equation}
(cf. \cite{ca,ma,me,mao}) in a different context.  Here one has a wonderful
theory developed primarily by S. Majid (cf. \cite{ma} and references there). 
Basically the idea here is that every quasitriangular Hopf algebra H has a braided
group analogue $\ul{H}$ of enveloping algebra type, given by the same algebra, unit,
and counit as H, but with a covariantized coproduct $\ul{\gD}$ which is braided
cocommutative.  The procedure $H\to\ul{H}$ is called transmutation.  On the other
hand any braided Hopf algebra B, in the braided category of left H modules for a 
quasitriangular H, gives rise to an ordinary Hopf algebra $B\lcrb\, H$ called the
bosonization of B; the modules of B in the braided category correspond to ordinary
modules of $B\lcrb\, H$.  In view of some theorems in \cite{ma} about bosonization
of the braided line, $U_q(s\ell(2)),$ etc. one expects perhaps to find relations to
other forms of transmuation as intertwining (cf. \cite{ca} for more on this).
In any event we recall that the R-matrix, arriving historically in quantum inverse
scattering theory, gives rise to quasitriangular Hopf algebra (i.e. quantum groups)
and thence to braiding.  Braiding, with its connections to knots, etc., seems to be
what is it's all about in connecting classical and quantum mathematics.
\\[3mm]\indent
The idea of intertwining operators plays a classically important role in the theory
of group representations and in particular there are connections between
intertwining operators, Hirota bilinear equations, and tau functions, which have
some noncommutative and quantum versions (see e.g.
\cite{gb,gc,ku,kz,kww,mbb,mbc,mbd,mbe,va,za}).  The idea here is that a
generic tau function can be defined as a generator of all the matrix elements of
$g\in G$ in a given highest weight representation of a universal enveloping 
type algebra.
The classical KP and Toda tau functions arise when G is a Kac-Moody algebra of level
$k=1$ (cf. also \cite{cx}).  In the case of quantum groups such ``tau" functions are
not number valued but take their values in noncommutative algebras (of functions on
the quantum group G).  The generic tau functions also satisfy bilinear Hirota type
q-difference equations, arising from manipulations with intertwining operators.  We
will not try to deal with this fascinating topic here but will pick it up in
\cite{ca}.  It does provide a source of natural q-difference equations related to
intertwining in a  quantum group context (there are also connections to bosonization
ideas in a  quantum mechanical spirit).
\\[3mm]\indent
Finally we will sketch some ideas about q-calculus following Wess and Zumino
\cite{wd,wzz} for example (there is an extensive literature on differential calculi
going back at least to Manin and Woronowicz (cf. \cite{ka,maa,wb}) and we do not try
to review this here (see however \cite{ca} for some of this); there seem to be
too many differential calculi and the Wess-Zumino approach develops matters directly
with the action of quantum derivatives acting on a quantum plane and direct
connections to the Heisenberg algebra.  This should set the stage for beginning a
construction of quantum transmutation following the classical patterns using partial
differential equations and spectral couplings of the author.  Further development is
planned for
\cite{ca}.
\\[3mm]\indent
Let us make a few remarks about discretization 
referring to \cite{ca,ch,cb,dae,hg,kad} for some background.
It has long been known that standard
calculus of the continuum is not accurate or even reasonable at the quantum level
(see e.g.
\cite{ch,cb}).  The language of quantum groups for example provides natural
discretizations, some of which are discussed in \cite{bi,bm,ca,ch,fw,
nb}, in connection with standard discretizations, and generally
noncommutativity and discrete physics have a natural compatibility (cf. \cite
{kad}).  It seems to me however
that there should be a natural discretization imposed via the ``principle" that time
is a secondary conception generated by the transfer of energy into information (cf.
\cite{ch,hg,mb}).  It also seems to be weird that time has been so successfully
employed as a dimension as in Minkowski space for example when it clearly has no
relation to dimensions at all (cf. \cite{ch,czb} 
and references there for other remarks on time).  In any
event as a function of energy we could have
$t=t(E)$ in many forms; in particular, thinking of a clock for example, time should
often appear most naturally as a discrete quantity.  Continuous transfer of energy
into time as information does not seem to be excluded and may be macroscopially
available; however we conjecture that at the quantum level the transfer ``must be"
discrete.  Then we could take here
$\gD t=t(E+\gD E)-t(E)$ or for simplicity in argument say $\gD t=f(E)\gD E$.  Now
consider kinetic energy corresponding to say $E=(1/2)\dot{x}^2=(1/2)p^2$ and write
\bq\label{1z}
\gD E=\frac{1}{2}\left(\frac{\gD x}{\gD t}\right)^2=\frac{1}{2}\left(\frac{\gD x}
{\gD E}\right)^2\left(\frac{\gD E}{\gD t}\right)^2=\frac{1}{2}\left(\frac{\gD x}{\gD
E}\right)^2f^{-2}(E)
\end{equation}
If now $f$ is a constant $T$ so that $\gD t=T\gD E$ (uniform time measurement) then
\eqref{1z} becomes
\bq\label{2z}
\frac{\gD t}{T}=\frac{1}{2}\left(\frac{\gD x}{\gD t}\right)^2\left(\frac{\gD t}
{\gD E}\right)^2T^{-2}=\frac{1}{2}\left(\frac{\gD x}{\gD t}\right)^2\Rightarrow
\end{equation}
$$\Rightarrow (\gD x)^2=\frac{2}{T}(\gD t)^3\sim |\gD x|=\sqrt{\frac{2}{T}}
|\gD t|^{3/2}$$
Hence if one is discretizing some standard differential equation involving time with
grid spacings $|\gD x|=\gs$ and $|\gD t|=\tau$ it appears natural to work with 
${\bf (Z1)}\,\,\gs=(2/T)^{1/2}\tau^{3/2}$.  We note that there are many
discretizations corresponding to a given PDE for example so there should be a choice
mechanism (perhaps); 
the alternative is to derive quantum PDE directly from a quantum group
context for example (which is undoubtedly the best recourse).  This will
be discussed later.
\\[3mm]\indent
We do not know of any calculations in this spirit by numerical analysts but this is not
an area of personal expertise.  Note however that if one takes $\tau=1$ and $\gs=\ga$
then ${\bf (Z2)}\,\,\ga\sim (2/T)^{1/2}$ which gives a measurement of the energy
conversion.  Thus if some discretization process requires a certain value of $\ga$ for 
implementation then a significance for $\ga$ arises which would then reveal the value  
of $T$.  For example in \cite{bm,nb} one works 
separately with discrete space
(and discrete time) quantization of the Schr\"odinger equation $\mf{S}$, leading to
quantum Hopf symmetry algebras $U_{\gs}(\mf{S})$ and $U_{\tau}(\mf{S})$.  This
involves nonstandard quantum deformations of $s\ell(2{\bf R})$ (Jordan or
h-deformations) and some twist maps which relate the two pictures.  These algebras
also come as limiting objects in a joint discretization when $\gs$ or $\tau$ tend to
zero.  However a combined q-group picture for joint discretization is not yet clear.
\begin{example}
Note for $q=exp(h)$ one has (for $\gD t\sim ht$ depending on $t$)
\bq\label{3z}
f(qt)=f(e^ht)\sim f((1+h)t)=f(t+\gD t);\,\,f(q^2t)=f(e^{2h}t)\sim f(t+2\gD t)
\end{equation}
so this suggests a dependency $\gD t=f(E)\gD E$ for a corresponding standard
discretization $\gD f=f(t+\tau)-f(t)$ to represent (approximately) a q-derivative
numerator $f(qt)-f(t)$.  If $t(E)$ were to be differentiable (with still discrete
transfer of energy) one could write (approximately) $\gD t/\gD E\sim t'(E)$ to get
$\gD t\sim t'(E)\gD E$.  then if $\gD t=ht$ one would need $ht\sim t'(E)\gD E$ which
could be realized via say $\gD E=c$ and $ct'=ht$ or $t=exp[(h/c)E]$.
Actually it seems increasingly attractive to simply look at the $q$ or $h$
discretizations arising naturally in a quantum mechanical context
with e.g. q-derivatives
$\pp_qf(t)=[f(qt)-f(t)]/[qt-t]$ and $q=exp(h)$ as the basic discretization involving
$\gD t=(q-1)t=((exp(h)-1)t$ (cf. also \cite{sk,tb}).
\end{example}

\section{SOME QUANTUM TRANSMUTATIONS}
\renewcommand{\theequation}{2.\arabic{equation}}

In order to prepare matters for further discussion we record a few
formulas from \cite{kf} (cf. also \cite{daa})
to illustrate what has already been done in related
directions and to give an exposure to q-analysis (for more details on q-calculus see
e.g. \cite{ca,kb,kj,ki}).
We use the notation 
\bq\label{44b}
{}_r\phi_s(a_1,a_2,\cdots a_r;b_1,b_2,\cdots b_s;q;z)=
\end{equation}
$$=\sum_0^{\infty}\frac{(a_1;q)_n(a_2;q)_n\cdots (a_r;q)_n}{(q;q)_n(b_1;q)_n\cdots 
(b_s;q)_n}\left[(-1)^nq^{\frac{n(n-1)}{2}}\right]^{1+s-r}z^n$$
for general hypergeometric functions
and write little q-Jacobi functions via
\bq\label{12z}
\phi_{\gl}(x;a,b;q)={}_2\phi_1(a\gs,a/\gs;ab;q;-(bx/a));\,\,\gl=\frac{1}{2}
(\gs+\gs^{-1})
\end{equation}
One has a general hypergeometric q-difference operator ($T_qf(x)=f(qx)$)
\bq\label{13z}
L=L^{(a,b)}=a^2\left(1+\frac{1}{x}\right)(T_q-id)+
\left(1+\frac{aq}{bx}\right)(T_q^{-1}-id)
\end{equation}
and the little q-Jacobi function satisfies ${\bf (Z3)}\,\,L\phi_{\gl}(\cdot;a,b;q)=
(-1-a^2+2a\gl)\phi_{\gl}(\cdot;a,b;q)$.  It is useful also to note that the little 
q-Jacobi functions are eigenvalues for eigenvalue $\gl$ of the operator
\bq\label{14z}
\mf{L}^{(a,b)}=\frac{a}{2}\left(1+
\frac{1}{x}\right)T_q-\left(\frac{a}{2x}+\frac{q}{2bx}\right)id+\frac{1}{2a}\left(
1+\frac{aq}{bx}\right)T_q^{-1}
\end{equation}
(i.e. $\mf{L}^{(a,b)}=(1/2a)L^{(a,b)}+(1/2)(a+a^{-1})$).  For simplicity one can assume
$a,b>0,\,\,ab<1,\,\,y>0$ but the results hold for a more general range of parameters
(cf. \cite{kf}).  The operator L is an unbounded symmetric operator on the Hilbert
space $\mf{H}(a,b;y)$ of square integrable sequences $u=(u_k)$ for $k\in{\bf Z}$ 
weighted via
\bq\label{15z}
\sum_{-\infty}^{\infty}|u_k|^2(ab)^k\frac{(-byq^k/a;q)_{\infty}}{(-yq^k;q)_{\infty}}
\end{equation}
where the operator L is initially defined on sequences with finitely many nonzero
terms.  Note that \eqref{15z} can be written as a q-integral by associating to u a
function f of $yq^{{\bf Z}}$ by $f(yq^k)=u_k$ and setting
\bq\label{16z}
\int_0^{\infty (y)}f(x)d_qx=y\sum_{-\infty}^{\infty}f(yq^k)q^k
\end{equation}
Then for $a=q^{(1/2)(\ga+\gb+1)}$ and $b=q^{(1/2)(\ga-\gb+1)}$ the sum in \eqref{15z}
can be written as
\bq\label{17z}
y^{-\ga-1}\int_0^{\infty(y)}|f(x)|^2x^{\ga}\frac{(-xq^{-\gb};q)_{\infty}}{(-x;q)_{\infty}}d_qx
\,\,\,\,(\Re\ga>-1)
\end{equation}
(note $b/a=q^{-\gb}$ and $ab=q^{\ga+1}$).
The spectral analysis of L, or 
equivalently of $\mf{L}^{(a,b)}$, can be carried out and leads to the transform
($\hat{u}=\mf{F}_{a,b,y}u$)
$$\hat{u}(\gl)=\sum_{k=-\infty}^{\infty}u_k\phi_{\gl}(yq^k;a,b;q)(ab)^k\frac
{(-byq^k/a;q)_{\infty}}{(-yq^k;q)_{\infty}}=y^{-\ga-1}\int_0^{\infty(y)}f(x)x^{\ga}
\phi_{\gl}(x;a,b;q)d_qx;$$
\bq\label{18z}
u_k=\int_{{\bf R}}(\mf{F}_{a,b,y}u)(\gl)\phi_{\gl}(yq^k;a,b;q)d\nu(\gl;a,b;y;q)
\end{equation}
where the measure is obtained from the c-function for expansions in little q-Jacobi
functions (cf. \cite{kf,kj,kg} - the formulas in \cite{kg} seem to differ but for 
structural purposes it doesn't matter).  The goal in
\cite{kf} is to establish a number of links between two little q-Jacobi function
transforms for different parameter values $(a,b,y)$ and revolves around 
transmutation kernels P satisfying 
\bq\label{19z}
(\mf{F}_{c,d,y}[\gd_tu])(\mu)=\int_{{\bf R}}(\mf{F}_{c,d,y}u)(\gl)P_t(\gl,\mu)d\nu
(\gl;a,b;y;q)
\end{equation}
where $(\gd_tu)_k=t^ku_k$ for an extra parameter t ($P_t\sim$ Poisson kernel for
$(a,b)=(c,d)$).  Similarly one studies the
transmuation kernel for the inverse transform
\bq\label{20z}
(\mf{F}_{c,d,y}^{-1}f)_{\ell}=\sum_{k=-\infty}^{\infty}(\mf{F}^{-1}_{c,d,y}f)_k
P_{k,\ell}(ab)^k\frac{(-byq^k/a;q)_{\infty}}{(-yq^k,q)_{\infty}}
\end{equation}
and one has the result
\bq\label{25z}
P_{k,\ell}(a,b,y;r,s)=\mf{F}^{-1}_{a,b,y}[\gl\to\phi_{\gl}(yq^{\ell}/s;ar,bs;q)]_k=
\end{equation}
$$=\int_{{\bf
R}}\phi_{\gl}(yq^{\ell}/s;ar,bs;q)\phi_{\gl}(yq^k;a,b;q)d\nu(\gl;a,b;y;q)$$
Another result involves the second order q-difference operator $\mf{L}^{(a,b)}$ as in
\eqref{14z} and this uses the Abel transform (cf. \cite{cw,kzz}) which has some nice
transmutation properties (cf. Section 3).  
Thus let $a,b\in{\bf C}/\{0\}$ and $\nu,\mu\in{\bf C}$ with
$|q^{\nu-\mu}b/a|<1$; then define the operator
\bq\label{28z}
(W_{\nu,\mu}(a,b)f)(x)=\frac{(-x;q)_{\infty}}{(-xq^{-\mu};q)_{\infty}}q^{-\mu^2}
\left(\frac{b}{a}\right)^{\mu}x^{\mu+\nu}\times
\end{equation}
$$\times\sum_{p=0}^{\infty}f(xq^{-\mu-p})q^{-p\nu}\frac{(q^{\nu};q)_{\infty}}
{(q;q)_{\infty}}{}_3\phi_2(q^{-p},q^{-\mu},-q^{1+\mu-\nu}a/bx;q^{1-p-\nu},-q^{\mu+1}/x;
q,q^{1-\mu}b/a)$$
for any function f with $|f(xq^{-p})|=O(q^{p(\gep+\nu)})$ for some $\gep>0$.  Then
$W_{\nu,\mu}(a,b)\ci\mf{L}^{(a,b)}=\mf{L}^{(aq^{-\nu},bq^{-\mu})}\ci W_{\nu,\mu}(a,b)$
on the space of compactly supported functions and for $\Phi_{\gs}$ the asympototically
free solution of $L\Phi_{\gs}(\cdot;a,b;q)=(-1-a^2+2a\gl)\Phi_{\gs}(\cdot;a,b;q)$ on
$yq^{{\bf Z}}$ ($\gl=(1/2)(\gs+\gs^{-1})$) one has
\bq\label{29z}
(W_{\nu,\mu}(a,b)\Phi_{\gs}(\cdot;a,b;q))(yq^k)=
y^{\mu+\nu}\frac{(a\gs,b\gs;q)_{\infty}}
{(aq^{-\nu}\gs,bq^{-\mu}\gs;q)_{\infty}}\Phi_{\gs}(yq^k;aq^{-\nu},bq^{-\mu};q)
\end{equation}
(note $(a,b;q)_{\infty}=(a;q)_{\infty}(b;q)_{\infty}$).
Further for $a,b>0,\,\,ab<1,\,\,\nu>0,$ and $\mu\in {\bf C}/{\bf Z}_{\leq 0}$ one
defines the operator
\bq\label{30z}
(A_{\nu,\mu}(a,b)f)(x)=\frac{(-bxq^{\mu}/a;q)_{\infty}}
{(-bxq^{\mu-\nu}/a;q)_{\infty}}\times
\end{equation}
$$\times\sum_{k=0}^{\infty}f(xq^{\mu+k})(ab)^k\frac{(q^{\nu},-xq^{\mu};q)_k}
{(q,-bxq^{\mu}/a;q)_k}{}_3\phi_2(q^{-k},q^{\mu},-bxq^{\mu-\nu}/a;q^{1-\nu-k},
-xq^{\mu};q,q)$$
for any bounded function f.  Then $\mf{L}^{(aq^{\nu},bq^{\mu})}\ci
A_{\nu,\mu}(a,b)=A_{\nu,\mu}(a,b)\ci\mf{L}^{(a,b)}$ on compactly supported functions
and
\bq\label{31z}
(A_{\nu,\mu}(a,b)\phi_{\gl}(\cdot;a,b;q))(x)=\frac{abq^{\nu+\mu};q)_{\infty}}
{(ab;q)_{\infty}}\phi_{\gl}(x;aq^{\nu},bq^{\mu};q)
\end{equation}
The operators $W_{\nu,\mu}(a,b)$ and $A_{\nu,\mu}(a,b)$ are referred to as q-analogues of
the generalized Abel transform.
One uses also the 
operator $W_{\nu}$ for $\nu\in{\bf C}$ acting on functions over
$[0,\infty)$ via
\bq\label{33z}
(W_{\nu}f)(x)=x^P{\nu}\sum_{\ell=0}^{\infty}f(xq^{-\ell})q^{-\ell\nu}
\frac{(q^{\nu};q)_{\ell}}{(q;q)_{\ell}}
\end{equation}
for $x\in [0,\infty)$, where one assumes that the infinite sum is absolutely
convergent if $\nu\not\in -{\bf Z}_{\geq 0}$.  For this one wants f sufficiently
decreasing on a q-grid tending to infinity, e.g. $f(xq^{-\ell})=O(q^{\ell(\nu+\gep)}$
for some $\gep>0$.  Note that for $\nu\in {\bf Z}_{\leq  0}$ the sum in \eqref{33z}
is finite and $W_0=id$ with $W_{-1}=B_q=(1/x)(1-T_q^{-1})$.  
This operator $W_{\nu}$ is a q-analogue of the Weyl
fractional integral operator for the Abel transform (cf. Section 3).  
Using the notation ${\bf (Z6)}
\,\,\int_a^{\infty}f(t)d_qt=a\sum_0^{\infty}f(xq^{-k})q^{-k}$ for the q-integral one
sees that for $n\in{\bf N}$ the operator $W_n$ is an iterated q-integral
\bq\label{34z}
(W_nf)(x)=\int_x^{\infty}\int_{x_1}^{\infty}\cdots\int_{x_{n-1}}^{\infty}f(x_n)
d_qx_nd_qx_{n-1}\cdots d_qx_1
\end{equation}
There are many other interesting formulas and results in \cite{kf} which we omit here;
note however
for $\mf{F}_{\rho}=\{f:\,[0,\infty)\to{\bf
C};\,|f(xq^{-\ell})|=O(q^{\ell\rho});\,\ell\to\infty;\,\forall x\in (q,1]\}$ one has
\begin{itemize}
\item
$W_{\nu}$ preserves the space of compactly supported functions
\item
$W_{\nu}:\,\mf{F}_{\rho}\to\mf{F}_{\rho-\Re\nu}$ for $\rho >\Re\nu>0$
\item
$W_{\nu}\ci W_{\mu}=W_{\nu+\mu}$ on $\mf{F}_{\rho}$ for $\rho>Re(\mu+\nu)>0$
\item
$W_{\nu}\ci B_q=B_q\ci W_{\nu}=W_{\nu-1}$ on $\mf{F}_{\rho}$ for $\rho>\Re \nu-1>0$
and 
$B_q^n\ci W_n=id$ for $n\in{\bf N}$ on $\mf{F}_{\rho}$ for $\rho>n$
\item
$\mf{L}^{(aq^{-\nu},b)}\ci W_{\nu}=W_{\nu}\ci\mf{L}^{(a,b))}$ for compactly supported
functions
\end{itemize}

\section{CLASSICAL TRANSMUTATION}
\renewcommand{\theequation}{3.\arabic{equation}}
\setcounter{equation}{0}

There are two classical approaches to transmutation of operators which should have
counterparts.  First (cf. \cite{cw,cv,cu,cx,ct,de,le} for details, hypotheses, and
examples), taking a special situation, if one has unique solutions to the partial
differential equation (PDE)) ${\bf (Z7)}\,\,P(D_x)\phi=Q(D_y)\phi$ with $\phi(x,0)=f(x)$
and e.g. $D_y\phi(x,0)=0$ (for suitable f) then we can define 
\bq\label{799}
Bf(y)=\phi(0,y);\,\,QBf=BPf
\end{equation}
(note that the boundary conditions can be generalized considerably).
To see this set $\psi=P(D_x)\phi$ so $[P(D_x)-Q(D_y)]\psi=0$ and $\psi(x,0)=Pf$ with
$\psi_y(x,0)=0$ so one can say that $\psi(0,y)=BPf$ while
$\psi(0,y)=P(D_x)\phi(x,y)|_{x=0}=Q(D_y)\phi(x,y)|_{x=0}=Q\phi(0,y)=QBf$.  There should be
some version of this for quantum operators.
\\[3mm]\indent
Secondly (using the above background) we got a lot of mileage from explicitly 
constructing
transmutation operators via eigenfunctions of P and Q with a spectral pairing 
(cf. \cite
{bv,cw,cv}).  Suppose for example that P and Q both have the same continuous 
spectrum $\gL$
with
\bq\label{800}
P(D_x)\phi_{\gl}^P(x)=\gl\phi_{\gl}^P;\,\,Q(D_x)\phi_{\gl}^Q=\gl\phi_{\gl}^Q
\end{equation} 
where e.g. ${\bf (Z8)}\,\,\phi_{\gl}^P(0)=\phi_{\gl}^Q(0)=1$ and
$D\phi_{\gl}^P(0)=D\phi_{\gl}^Q(0)=0$ (again this can be generalized).  Let the related
spectral measures be $d\mu_P$ and $d\mu_Q$ with Fourier type recovery formulas (for
suitable f)
\bq\label{801}
\hat{f}_P(\gl)=\int
f(x)\gO^P_{\gl}(x)dx;\,\,f(x)=\int_{\gL}\hat{f}_P(\gl)\phi^P_{\gl}(x)d\mu_P;
\end{equation}
$$\hat{f}_Q(\gl)=\int
f(x)\gO^Q_{\gl}(x)dx;\,\,f(x)=\int_{\gL}\hat{f}_Q(\gl)\phi_{\gl}^Q(x)d\mu_Q$$
Here for second order differential operators ${\bf (Z9)}\,\,(Au')'/A+ qu=\gl u$ (q real),
arising in the treatment of many special functions, one has generically
$\gO_{\gl}^P(x)=A\phi_{\gl}^P(x)$ and $d\mu_P$ is usually absolutely continuous.
We note that for suitable $f,\,g$ and P as in {\bf (Z9)} one has formally ${\bf
(Z10)}\,\,\int[P(D)f]gdx=\int[(Af')'/A+qf]g=\int [-Af'(g/A)'+qfg]=\int f\{[(g/A)'A]'
+qg\}$ so
$P^*(D)g\sim [A(g/A)']'+qg$.  In particular we have ${\bf (Z11)}\,\,
P^*(D)\gO_{\gl}^P=\gl\gO^P_{\gl}$.  Now one can define
\bq\label{803}
\phi(x,y)=\int\phi_{\gl}^P(x)\phi_{\gl}^Q(y)\hat{f}_P(\gl)d\mu_P
\end{equation}
and check that formally $P(D_x)\phi=Q(D_y)\phi$ with $\phi(x,0)=\int\phi_{\gl}^P(x)\hat
{f}_P(\gl)d\mu_P=f(x)$ and $\phi_y(x,0)=0$.  Consequently as in \eqref{799} one can write
${\bf (Z12)}\,\,\phi(0,y)=Bf(y)=\int\phi_{\gl}^Q(y)\hat{f}_P(\gl)d\mu_P$.
\\[3mm]\indent
We write down now a few more general features for completeness and in order to exhibit
some group (and perhaps quantum group) theoretic content.  In particular one can define a
generalized translation operator $T_{\xi}^xf(\xi)=U(x,\xi)$ via
\bq\label{804}
U(x,\xi)=\int_{\gL}\hat{f}_P(\gl)\phi^P_{\gl}(x)\phi^P_{\gl}(\xi)d\mu_P
\end{equation}
so $P(D_x)U=P(D_{\xi})U$ with $U(x,\xi)=U(\xi,x)$ and then one sets
($Bf(y)=<\gb(y,x),f(x)>$)
\bq\label{805}
\gb(y,x)=<\gO^P_{\gl}(x),\phi_{\gl}^Q(y)>_P=\int\gO^P_{\gl}(x)\phi_{\gl}^Q(y)d\mu_P;
\,\,\phi(x,y)=<\gb(y,\xi),U(x,\xi)>
\end{equation}
To see that this is formally the same as \eqref{803} we write from \eqref{801}
the relations ${\bf
(Z13)}\,\,\hat{f}_P(\gl)=\int\gO^P_{\gl}(x)\left(\int_{\gL}\hat{f}_P(\gz)
\phi_{\gz}^P(x)d\mu_P(\gz)\right)dx\sim\int_{\gL}\hat{f}_P(\gz)d\mu_P(\gz)
(\int\gO^P_{\gl}(x)\phi_{\gz}^P(x)dx)$ which implies that
$d\mu_P(\gz)\int\gO^P_{\gl}(x)\phi^P_{\gz}(x)dx\sim
\gd(\gl-\gz)d\gz$ (Darboux-Christoffel arguments can also be used here) while
${\bf (Z14)}\,\,f(x)=\int_{\gL}\phi^P_{\gl}(x)\left(\int
f(\xi)\gO^P_{\gl}d\xi\right)d\mu_P (\gl)\sim\int
f(\xi)\left(\int\phi_{\gl}^P(x)\gO^P_{\gl}(\xi)d\mu_P(\gl)\right)d\xi$ implies
$\int\phi_{\gl}^P(x)\gO^P_{\gl}(\xi)d\mu_P(\gl)\sim\gd(x-\xi)$.  Then
\eqref{805} becomes
\bq\label{806}
<\gb(y,\xi),U(x,\xi)>\sim\int\int\int\gO^P_{\gl}(x)\phi^Q_{\gl}(y)
d\mu_P(\gl)\hat{f}(\gz)\phi^P_{\gz}(x)\phi^P_{\gz}(\xi)d\mu(\gz)d\xi\sim
\end{equation}
$$\sim\int\int \phi_{\gz}^P(x)\phi_{\gl}^Q(y)\hat{f}(\gz)d\mu_P(\gl)d\mu_P(\gz)
\int\gO^P_{\gl}(\xi)\phi_{\gz}^P)(\xi)d\xi\sim$$
$$\sim\int\int\phi_{\gz}^P(x)\phi^Q_{\gl}(y)\hat{f}(\gz)\gd(\gz-\gl)d\mu_P(\gl)=
\int\phi^P_{\gl}(x)\phi^Q_{\gl}(y)\hat{f}(\gl)d\mu_P(\gl)$$
Another useful observation from \eqref{805} is (recall from {\bf (Z11)}
$P^*(D)\gO^P_{\gl}=\gl\gO^P_{\gl}$ implies that $P^*(D_{\xi})\gb(y,\xi)=P(D_y)\gb(y,\xi)$) 
\bq\label{807}
P(D_x)\phi=<\gb(y,\xi),P(D_x)U(x,\xi)>=
\end{equation}
$$=<\gb(y,\xi),P(D_{\xi})U(x,\xi)>=<P^*(D_{\xi})\gb(y,\xi),U(x,\xi)>=P(D_y)\phi$$
with $\phi(x,0)=<\gb(0,\xi),U(x,\xi)>$ and $\phi(0,y)=<\gb(y,\xi),U(0,\xi)>$.  But from
{\bf (Z15)}, $\gb(0,\xi)=<\gO^P_{\gl}(\xi),1>=\int\gO^P_{\gl}(\xi)d\mu_P(\gl)=\gd(\xi)$
and from
\eqref{804} and \eqref{801} it follows that
$U(0,\xi)=\int_{\gL}\hat{f}_P(\gl)\phi^P_{\gl}(\xi)d\mu_P(\gl)=
f(\xi)$ so $\phi(x,0)=U(x,0)=f(x)$ and $\phi(0,y)=<\gb(y,\xi),f(\xi)>=Bf(y)$ as desired.
\\[3mm]\indent
Now for the ``canonical" development of the author in \cite{cw} we  
write e.g. $A(x)\sim \gD_P(x)$ in $P(D)$ of {\bf (Z9)} and set
\bq\label{50z}
P(D)u=\frac{(\gD_Pu')'}{\gD_P}+pu;\,\,Q(D)u=\frac{(\gD_Qu')'}{\gD_Q}+qu
\end{equation}
and from \eqref{800}-\eqref{801} we repeat (with generic conditions $\phi^P_{\gl}(0)
=1$ and $D\phi^P_{\gl}(0)=0$)
\bq\label{51z}
P(D)\phi^P_{\gl}(x)=\gl\phi_{\gl}^P(x);\,\,P^*(D)\gO_{\gl}^P(x)=\gl\gO^P_{\gl};
\,\,P^*(D)f=\left[\gD_P\left(\frac{f}{\gD_P}\right)'\right]'+pf
\end{equation}
where $\gO^P_{\gl}(x)=\gD_P(x)\phi_{\gl}^P(x)$.  The following transforms are then
relevant where we write $<f,\phi>\sim\int f(x)\phi(x)dx$ over some range and
$<F,\psi>_{\nu}\sim\int F(\gl)\psi(\gl)d\nu(\gl)$ over some range (the ranges may
be discrete or contain discrete sections).  We also note that many measure pairings
involving $d\nu(\gl)$ are better expressed via distribution pairings with a
generalized spectral function (distribution) as in \cite{cw,mba} (note also that 
transmutation does not require that P and Q have the same spectrum - cf. 
\cite{bv,cw}).  Thus we write
for suitable $f,F$ (and $d\go$ the spectral measure associated with Q)
\begin{enumerate}
\item[{\bf A}.]
$\mf{P}f(\gl)= <f(x),\gO^P_{\gl}(x)>;\,\,\,\mf{Q}f(\gl)=<f(x),\gO_{\gl}^Q(x)>$
\item[{\bf B}.]
${\mc P}f(\gl)=<f(x),\phi_{\gl}^P(x)>;\,\,\,{\mc Q}f(\gl)=<f(x),\phi_{\gl}^Q(x)>$
\item[{\bf C}.]
$\tl{\mf{P}}F(x)=<F(\gl),\phi_{\gl}^P>_{\nu};\,\,\,\tl{\mf{Q}}F(x)=<F(\gl),
\phi_{\gl}^Q(x)>_{\go};\,\,\tl{\mf{P}}\sim\mf{P}^{-1};\,\,\tl{\mf{Q}}\sim\mf{Q}^{-1}$
\item[{\bf D}.]
$\tl{{\mc P}}F(x)=<F(\gl),\gO_{\gl}^P(x)>_{\nu};\,\,\,\tl{{\mc Q}}F(x)=<F(\gl),
\gO_{\gl}^Q>_{\go};\,\,\tl{{\mc P}}\sim {\mc P}^{-1};\,\,\tl{{\mc Q}}\sim {\mc Q}^{-1}$
\item[{\bf E}.]
${\bf P}F(x)=<F(\gl),\phi_{\gl}^P(x)>_{\go};\,\,\,{\bf Q}F(x)=<F(\gl),\phi^Q_{\gl}
(x)>_{\nu}$
\end{enumerate}
One notes the mixing of eigenfunctions and measures in {\bf E} and this gives rise to
the transmutation kernels (${\mc B}=B^{-1}$)
\bq\label{52z}
ker(B)=\gb(y,x)=<\gO^P_{\gl}(x),\phi^Q_{\gl}(y)>_{\nu};\,\,\,ker({\mc
B})=\gag(x,y)=<\phi_{\gl}^P(x),\gO^Q_{\gl}(y)>_{\go}
\end{equation}
leading to 
\bq\label{53z}
B={\bf Q}\ci\mf{P};\,\,\,{\mc B}={\bf P}\ci \mf{Q}
\end{equation}
\indent
This is a very neat picture and it applies in a large number of interesting situations
(cf. \cite{cw,cv,cu,cx}).  We indicate some further features now comprising a classical
development on which the example of \cite{kf} 
(discussed briefly in Section 2), and much more, can be
predicated. The proofs in the classical situation involve a heavy dose of Paley-Wiener
theory and Fourier ideas.  Thus from \cite{cw}, Example 9.3, we take ${\bf
(Z16)}\,\,\gD_Q=\gD_{\ga,\gb}=(e^x-e^{-x})^{2\ga+1}(e^x+e^{-x})^{2\gb +1}$ with 
$\rho=\ga+\gb+1$.  Then for $\ga\ne -1,-2,\cdots$ one has
\bq\label{54z}
\phi_{\gl}^Q(x)=\phi^{\ga,\gb}_{\gl}(x)=F((1/2)(\rho+i\gl),(1/2)(\rho-i\gl),\ga+1,
-sh^2x)
\end{equation}
where $sh\sim sinh$ and $F\sim {}_2F_1$ is the standard hypergeometric function,
extended to ${}_2\phi_1$ in \eqref{44b} for the q-theory.
The related ``Jost" functions are
\bq\label{55z}
\Phi_{\gl}^Q(x)=(e^x-e^{-x})^{i\gl-\rho}F\left(\frac{1}{2}(\gb-\ga+1-i\gl),
\frac{1}{2}(\gb+\ga+1-i\gl),1-i\gl,-sh^{-2}x\right)
\end{equation}
for $\gl\ne -i,-2i,\cdots$.  Here $\Phi_{\gl}^Q\sim exp(i\gl-\rho)x$ as $x\to\infty$
and one can write ${\bf (Z17)}\,\,\phi_{\gl}^Q=c_Q\Phi_{\gl}^Q+\bar{c}_Q\Phi_{-\gl}^Q$
for $\gl\ne 0,\pm i,\pm 2i,\cdots$; here $c_Q$ is the Harish-Chandra c-function which
also generates the measure $d\go_Q\sim c|c_Q|^2d\gl$ (see \cite{kf} for the discrete
version of this).  For the Abel and Weyl transformations we
recall from \cite{cw} that the Fourier-Jacobi transform for $f\in C_0^{\infty}$ is
defined via
\bq\label{56z}
\hat{f}_{\ga,\gb}(\gl)=\frac{\sqrt{2}}{\gG(\ga+1)}\int_0^{\infty}f(t)\phi_{\gl}^
{\ga,\gb}(t)\gD_{\ga,\gb}(t)dt=\frac{\sqrt{2}}{\gG(\ga+1)}\mf{Q}f(\gl)
\end{equation}
Then, using some identities for hypergeometric functions one can show that
\bq\label{57z}
\hat{f}_{\ga,\gb}(\gl)=\frac{2}{\pi}\int_0^{\infty}F_{\ga,\gb}[f]Cos(\gl s)ds
\end{equation}
where $F_{\ga,\gb}[f]$ is the Abel transform ${\bf
(Z18)}\,\,F_{\ga,\gb}[f](x)=\int_0^{\infty}f(t)A(s,t)dt$ and $A(s,t)$ has an explicit
form in terms of a particular ${}_2F_1$.
\\[3mm]\indent
{\bf REMARK 3.1}.  Note that \eqref{57z} is a special case of a general formula in 
\cite{cw}.  It also has a version in the theory of Lie groups and symmetric spaces
where $exp(-\rho s)F_{\ga,\gb}[f](s)$ can be interpreted as a Radon transform of a
radial function $f$ and one can write in standard Lie theory notation
\bq\label{58z}
F_f(a)=e^{\rho(log a)}\int_Nf(an)dn;\,\,F^*(\gl)=\int_AF(a)a^{-i\gl(log a)}da;
\end{equation}
$$\tl{f}(\gl)=\int_Gf(x)\phi_{-\gl}(x)dx$$
Then $\tl{f}=(F_f)^*\sim \hat{f}_{\ga,\gb}(\gl)$ in \eqref{57z}.  The transmutation
version of this has the form ${\mc P}F_Q[f]=\mf{Q}f$.  Indeed following \cite{cw} we
can write for the Fourier-Jacobi situation ($\tl{c}_{\ga,\gb}=2\sqrt{\pi}c_{\ga,\gb}/
\gG(\ga+1)$)
\bq\label{61z}
\check{g}_{\ga,\gb}(t)=\frac{1}{\sqrt{2\pi}}\int_{-i\eta-\infty}^{i\eta+\infty}g(\gl)
\frac{\phi_{\gl}^{\ga,\gb}(t)}{\tl{c}_{\ga,\gb}(-\gl)}d\gl
\end{equation}
for suitable parameter values.  Then one exploits known relations for the Cosine
transform (corresponding to $(\ga,\gb)=(-1/2,-1/2)$) to arrive eventually at a formula
$(\check{g}_{\ga,\gb}\sim f,\,\,\hat{f}\sim g$)
\bq\label{62z}
\hat{f}(\gl)=\left(\frac{2}{\pi}\right)^{1/2}\int_0^{\infty}F_{\ga,\gb}[f](s)Cos(\gl s)
ds
\end{equation}
corresponding to ${\mc P}F_Q[f](\gl)=\mf{Q}f(\gl)$.  Transforms based on ${\bf (Z19)}
\,\,\Psi_{\gl}^Q(x)=\Phi_{\gl}^Q(x)/c_Q(-\gl)$ are important in the theory of the
Marchenko equation and are discussed below (cf. \cite{cw} for an extensive treatment).
\\[3mm]\indent
One can also use the Weyl fractional integral operators of the form ${\bf
(Z20)}\,\,W_{\mu}[g](y)=\gG(\mu)^{-1}\int_y^{\infty}g(x)(x-y)^{\mu-1}dx$ for
$\Re\mu>0$ and $g\in C_0^{\infty}$.  Then this can be extended as an
entire function in $\mu\in{\bf C}$ with ${\bf (Z21)}\,\,W_{\mu}\ci W_{\nu}=W_{\mu+\nu}$
where $W_{\mu}[g](y)\in C_0^{\infty},\,\,W_0=id,$ and $W_{-1}[g]=-g'$.
For $f\in C_0^{\infty}$ and $Re\nu>0$ one also defines now for the spherical function
situation of $\phi_{\gl}^{\ga,\gb}$ ($ch\sim cosh$)
\bq\label{59z}
\mf{W}_{\mu}^{\gs}[f](s)=\gG(\mu)^{-1}\int_s^{\infty}f(t)[ch(\gs t)-ch(\gs
s)]^{\mu-1}d(ch(\gs t)
\end{equation}
which can be extended to $\mu\in{\bf C}$ as an entire function with
$(\mf{W}_{\mu}^{\gs})^{-1}=\mf{W}_{-\mu}^{\gs}$.  Further one has (cf. \cite{cw})
\bq\label{60z}
F_{\ga,\gb}[f]=2^{3\ga+(3/2)}\mf{W}^1_{\ga-\gb}\ci \mf{W}^2_{\gb+(1/2)}[f];\,\,F^{-1}_
{\ga,\gb}=2^{-3\ga-(3/2)}\mf{W}^2_{-\gb-(1/2)}\ci\mf{W}^1_{\gb-\ga}
\end{equation}
\indent
We mention now formally a few additional formulas for suitable $f,g$.  If one writes 
$Bf(y)=<\gb(y,x),f(x)>$ (cf. \eqref{805}) and ${\mc B}g(x)=<\gag(x,y),g(y)>\,\,
({\mc B}=B^{-1})$ then setting ${\mc B}^*g(x)=<\gag(x,y),g(x)>$ 
and $B^*f(x)=<\gb(y,x),f(y)>$ there results
\bq\label{63z} 
{\mc P}B^*f={\mc Q}f;\,\,{\mc Q}{\mc B}^*g={\mc P}g
\end{equation}
(note ${\bf (Z22)}\,\,\gb(y,x)=<\gO^P_{\gl}(x),\phi^Q_{\gl}(y)>_{\nu}$ and $\gag(x,y)=
<\phi_{\gl}^P(x),\gO^Q_{\gl}(y)>_{\go}$).  There are also Parseval formulas of the form
\bq\label{64z}
<R,{\mc Q}f{\mc Q}g>_{\gl}=<\gD_Q^{-1/2}f,\gD_Q^{-1/2}g>;\,\,<R,\mf{Q}f\mf{Q}g>_{\gl}=
<\gD_Q^{1/2}f,\gD_Q^{1/2}g>
\end{equation}
where $R$ is a spectral function (distribution) associated to $d\go_Q(\gl)$ (cf.
\cite{cw} for details).  Evidently one has ${\bf (Z23)}\,\,B\phi_{\gl}^P=\phi_{\gl}^Q$
and from (cf. \eqref{53z}) ${\mc B}={\bf P}\ci\mf{Q}=B^{-1}=[{\bf
Q}\mf{P}]^{-1}=\tl{\mf{P}}{\bf Q}^{-1}$ we get
\bq\label{65z}
{\bf Q}^{-1}=\mf{P}{\bf P}\mf{Q};\,\,{\bf P}^{-1}=\mf{Q}{\bf Q}\mf{P}
\end{equation}
\indent
One goes next to the Gelfand-Levitan (GL) and Marchenko (M) equations which are of
fundamental importance in inverse scattering theory for example (cf. \cite{cw,cu,
cv,cx,cag,fu}).
One determines an adjoint operator $B^{\#}=\tl{{\mc B}}$ by the rule
\bq\label{66z}
<\gD_Q(y)v(y),Bu(y)>=<v(y),\gD_Q(y)<\gb(y,x)\gD_P^{-1}(x)\gD_P(x),u(x)>>=
\end{equation}
$$=<\gD_P(x)u(x),<\tl{\gag}(x,y),v(y)>>=<\gD_P(x)u(x),\tl{{\mc B}}v(x)>$$
where $\tl{\gag}=ker\tl{{\mc B}}$ is given via
\bq\label{67z}
\tl{\gag}(x,y)=\gD_Q(y)\gD_P^{-1}(x)\gb(y,x)=<\phi_{\gl}^P(x),\gO_{\gl}^Q(y)>_{\nu}
\end{equation}
We write also ${\mc B}^{\#}=\tl{B}$ where 
$$<\gD_P(x)u(x),{\mc B}v(x)>=<u,\gD_P<\gag(x,y)\gD_Q^{-1}\gD_Q,v>>=<\gD_Qv(y),
<\tl{\gb}(y,x),u(x)>>;$$
\bq\label{68z}
\tl{\gb}(y,x)=\gD_Q^{-1}(y)\gD_P(x)\gag(x,y)=<\gO_{\gl}^P(x),\phi_{\gl}^Q(y)>_{\go}
\end{equation}
One shows that $\tl{B}$ can also be characterized via a Cauchy problem as in Section 1.
Thus if $T^x_{\xi}$ is the generalized translation associated with P
(cf. \eqref{804} - i.e. $T^x_{\xi}f=
<\hat{f}_P(\gl),\phi_{\gl}^P(x)\phi_{\gl}^P(\xi)>_{\nu}$ with
$\hat{f}_P(\gl)=\mf{P}f(\gl)$) then set $\phi(x,y)=<\tl{\gb}(y,\xi),T^x_{\xi}f>$
and after some calculation one gets $P(D_x)\phi=Q(D_y)\phi,\,\,\phi(x,0)=\mf{A}f(x),$
and
$\phi_y(x,0)=0$ where ${\bf (Z24)}\,\,\mf{A}f(x)=<\mf{A}(x,\xi),f(\xi)>$ with
$\mf{A}(x,\xi)=<\gO^P_{\gl}(\xi),\phi_{\gl}^P(x)>_{\go}$ (since $\mf{A}$ commutes with
$P(D)$ one can apply the PDE method of Section 1 - cf. \cite{cw}).  One result of this
now is that ${\bf (Z25)}\,\,\tl{B}f(y)=<\phi_{\gl}^Q(y),\mf{P}f(\gl)>_{\go}=\tl{{\bf
Q}}\mf{P}f(y)$.  In a similar manner one can characterize $\tl{{\mc B}}$ via a Cauchy
problem as $\psi(x,0)=\tl{{\mc B}}g(x)$ where
$P(D_x)\psi=Q(D_y)\psi$ with $D_x\psi(0,y)=0$ and $\psi(0,y)=\check
{\mf{A}}g(y)$; here ${\bf
(Z26)}\,\,\check{\mf{A}}g(y)=<\check{\mf{A}}(y,\eta),g(\eta)>$ with $\check{\mf{A}}
(y,\eta)=<\gO^Q_{\gl}(\eta),\phi_{\gl}^Q(y)>_{\nu}$.  Using such operators one arrives
at a generalized G-L equation in the form
\bq\label{69z}
B\ci\mf{A}=\tl{B}\equiv\, <\gb(y,t),\mf{A}(t,x)>=\tl{\gb}(y,x)
\end{equation}
which can be reformulated as
\bq\label{70z}
\int<\gO^P_{\gl}(t),\phi^Q_{\gl}(y)>_{\nu}<\gO^P_{\gz}(x),\phi_{\gz}^P(t)>_{\go}dt=
<\gO^P_{\mu}(x),\phi_{\mu}(y)>_{\go}
\end{equation}
We recall also that the kernels satisfy triangularity theorems (arising classically
via Paley-Wiener type theorems) and in particular one proves that $\tl{\gb}(y,x)=0$
for $y>x$ and $\gb(y,x)=0$ for $x>y$.  We note also that for $d\nu=\hat{\nu}d\gl$ and
$d\go=\hat{\go}d\gl$ one has formally ${\bf
(Z27)}\,\,\tl{B}\phi_{\gl}^P=W\phi_{\gl}^Q$ where $W=\hat{\go}/\hat{\nu}$.  The
triangularity leads to the classical appearance of the GL theory in the form
$\gb(y,x)=\gd(x-y)+K(y,x)$ with $K(y,x)=0$ for $x>y$ and \eqref{69z} takes the well
known form for $x<y$
\bq\label{72z}
\tl{\gb}(y,x)=\mf{A}(y,x)+\int_0^yK(y,\xi)\mf{A}(\xi,x)d\xi=0
\end{equation}
There are many variations and examples given in \cite{cw} for instance (cf. also
\cite{cag} for 8 or 9 derivations of GL equations).
\\[3mm]\indent
The Marchenko equation is also used in inverse scattering theory and we have cast this
in a form also giving a general setting for Kontorovich-Lebedev (KL) theory (cf.
\cite{cw,cah,cai} for example).  The Marchenko equation is somewhat more complicated
from an operator point of view and we only supply here a few operator formulas (see
e.g. \cite{cag} for a more relaxed presentation with a more classical touch).
First (cf. {\bf (Z27)}) one notes that for
$d\nu_P\sim\hat{\nu}_Pd\gl$ and $d\go_Q=\hat{\go}_Qd\gl$ one can write for
$\Psi_{\gl}^{P,Q}=\Phi_{\gl}^{P,Q}/c_{P,Q}(-\gl)$ the formulas
\bq\label{72z}
\gb(y,x)=\frac{\gD_P(x)}{2\pi}\int_{-\infty}^{\infty}
\Psi_{\gl}^P(x)\phi^Q_{\gl}(y)d\gl;\,\,\tl{\gb}(y,x)=\frac{\gD_P(x)}{2\pi}
\int_{-\infty}^{\infty}\Psi_{\gl}^Q(y)\phi_{\gl}^P(x)d\gl
\end{equation}
(indicating nicely the role reversal between $x$ and $y$).  In particular one can show
that under suitable hypotheses ${\bf
(Z28)}\,\,\tl{B}[\Psi_{\gl}^P(x)](y)=\Psi^Q_{\gl}(y)$ and
$\tl{\gb}(y,x)=(1/2\pi)\int_{-\infty}^{\infty}\gO_{\gl}^P(x)\Psi_{\gl}^Q(y)d\gl$ as
well as (cf. \eqref{63z}) ${\bf (Z29)}\,\,{\mc B}^*[\gD_Pf]=\gD_Q\tl{B}f$ and
$B^*[\gD_Qf]=\gD_P\tl{{\mc B}}f$.  For certain $P,\,Q$ one can also generate
transmutations (acting on suitable functions) via kernels of the form 
${\bf (Z30)}\,\,\check{\gb}(y,x)=(1/2\pi)\int_{-\infty}^{\infty}\Phi_{\gl}^Q(y)
(\phi_{\gl}^P(x)/c_P(-\gl))d\gl$ which can be put in the form
${\bf
(Z31)}\,\,\check{\gb}(y,x)=(1/2\pi)\int_{-\infty}^{\infty}
\Phi_{\gl}^Q(\gl)\rho(\gl)\Sigma_{\gl}^P(x)d\gl$ for a full line theory.
Here we note (without discussion) that for suitable operators ${\bf
(Z32)}\,\,(c_P(\gl)/c_P(-\gl))\Phi_{\gl}^P(x)+\Phi_{-\gl}^P(x)=
\rho(\gl)[\Sigma_{\gl}^P(x)+\Phi_{\gl}^P(x)]$
while $\Phi_{-\gl}^P=\rho\Sigma_{\gl}^P+A\Phi_{\gl}^P$ and
$\Sigma_{-\gl}^P=\rho\Phi_{\gl}^P+A\Sigma_{\gl}^P$ (see \cite{cw} for details and
discussion).  The KL type inversion has the form $(\hat{f}(\gl)=\mf{Q}f(\gl)$ here)
\bq\label{73z}
f(x)=\frac{1}{2\pi}\int_{-\infty}^{\infty}\hat{f}(\gl)\Psi^Q_{\gl}(x)d\gl
\end{equation}
There is more discussion of KL type inversion in \cite{cw,cag,cah,cai}).  We note also
that an inverse $\check{{\mc B}}$ to $\check{B}$ is obtained via the kernel
\bq\label{74z}
\check{\gag}(x,y)=\frac{1}{2\pi}\int_{-\infty}^{\infty}
\Phi_{\gl}^P(x)\Psi^Q_{\gl}(y)d\gl
\end{equation}
\indent
Now the kernel for the Marchenko equation will involve terms
($s_Q(\gl)=c_Q(\gl)/c_Q(-\gl)$)
\bq\label{75z}
S(t,x)=\frac{1}{2\pi}\int_{-\infty}^{\infty}s_Q(\gl)\Phi_{\gl}^P(t)\Phi_{\gl}^P(x)d\gl;
\end{equation}
$$J(t,x)=\gd(t-x)+\frac{1}{2\pi}\int_{-\infty}^{\infty}A(\gl)
\Phi_{\gl}^P(t)\Phi_{\gl}^P(x)d\gl=
\gd(x-t)+M(t,x)$$
Here ${\bf (Z33)}\,\,J(t,x)=(1/2\pi)\int_{-\infty}^{\infty}\Phi_{\gl}^P(t)\Phi_{-\gl}^P
(x)d\gl$.  Now one can define
\bq\label{76z}
\hat{\gb}(y,\xi)=\frac{1}{2\pi}\int_{-\infty}^{\infty}\frac{1}{c_Q(-\gl)}\Phi_{\gl}^P
(\xi)\Phi_{\gl}^Q(y)d\gl
\end{equation}
and $\hat{\gb}(y,\xi)=0$ for $\xi>y$.
Here $\hat{B}=B{\mc H}^*$ where 
\bq\label{77z}
{\mc
H}(x,s)=\frac{1}{2\pi}\int_{-\infty}^{\infty}\frac{c_P(-\gl)}{c_Q(-\gl)}
\Sigma_{\gl}^P(s)\Phi_{\gl}^P(x)\rho d\gl
\end{equation}
\indent
The ingredients are now in place but not motivated (see \cite{cw,cag,cah,cai} for
motivation, examples, and further development).  Abstractly now the idea of the 
Marchenko equation is to relate B and $\tl{B}$ via $\check{B}$.  One obtains
first $\tl{B}=\check{B}{\mc H}$ and then the generalized extended Marchenko equation 
is ${\bf (Z34)}\,\,\hat{B}=
B{\mc H}^*=\tl{B}\tl{\mf{A}}{\mc H}^*=
\check{B}({\mc H}\tl{\mf{A}}{\mc H}^*)$
where from \eqref{69z} $\tl{B}=B\mf{A}\Rightarrow
B=\tl{B}\mf{A}^{-1}=\tl{B}\tl{\mf{A}}$ ($\tl{\mf{A}}=\mf{A}^{-1}$ and ${\mc
H}^*(x,s)={\mc H}(s,x)$). Note also ${\bf
(Z35)}\,\,\tl{\mf{A}}(t,x)=<\phi_{\gl}^P(x),\phi_{\gl}^P(t)W^{-1}(\gl)>_{\nu}$ where
$W^{-1}(\gl)=\hat{\nu}(\gl)/\hat{\go}(\gl)$.  Further
${\bf (Z36)}\,\,\check{B}^{-1}\hat{B}\sim {\mc H}\tl{\mf{A}}{\mc H}^*$ has kernel
$S+J$ and we note that the classical scattering theory
result emerges from $\hat{\gb}(y,x)=0$ for $x>y$
in the form
\bq\label{78z}
0=<\hat{\gb}(y,t),S(t,x)+J(t,x)>=\check{\gb}(y,x)+
\int_y^{\infty}\check{\gb}(y,t)[S(t,x)+M(t,x)]dt
\end{equation}
(cf. \cite{cw,cag} for more on relations to scattering theory).
One can discuss here also upper-lower operator factorizations and a parallel structure
for GL and Marchenko equations (cf. \cite{cw} and references there).  In fact for 
suitable operators one can think of $W(\gl)=\hat{\go}(\gl)/\hat{\nu}(\gl)$ as GL data
and $W^{-1}(\gl)$ as Marchenko data.
\\[3mm]\indent
We conjecture that much of this machinery can be rephrased in a q-analysis form and 
in this direction we
go next to some basic material concerning the quantum plane, $SL_q(2)$ and/or
$GL_q(2)$, plus differential calculus \`a la Wess-Zumino.  Eventually one then would
look for a completely
q-group formulation and connections to quantum tau functions etc. as mentioned in
Section 1. We will only make preliminary remarks toward implementation in this paper
but will try to give enough of the appropriate quantum group framework to make further
development natural.

\section{SOME Q-CALCULUS}
\renewcommand{\theequation}{4.\arabic{equation}}
\setcounter{equation}{0}

For our purposes the most attractive approach to q-calculus (following
\cite{cqc,cqd,oh,sz,wd,wzz}) starts with a quantum object, the Heisenberg algebra, and
builds a differential calculus around it, subsequently constructing the associated
quantum groups, quantum planes, integrals, etc.  We will sketch some of this here
following the superb exposition in \cite{wzz} to which we refer for details and a
more thorough treatment (cf. also \cite{ca}).
\subsection{{\bf The Heisenberg algebra.}}
The q-deformed Heisenberg algebra $\mf{H}_q$ has relations
\bq\label{1a}
q^{1/2}xp-q^{-1/2}px=i\gL;\,\,\gL p=qp\gL;\,\,\gL x=q^{-1}x\gL;\,\,q\ne 1;\,\,q\in
{\bf R}
\end{equation}
where $x$ and $p$ (position and momentum) are selfadjoint as operators (giving real
eigenvalues and a complete set of eigenfunctions). Thus we want a $*$-algebra and an
antilinear involution $a\to\bar{a}$ in the algebra (corresponding to $*$ in an
operator representation).  For this one needs to extend the algebra \eqref{1a}
by conjugate elements $\bar{x}=x,\,\,\bar{p}=p,\,\,\bar{q}=q$, and $\bar{\gL}=\gL^{-1}$
($\gL$ is unitary) plus $x^{-1}$.  Then ordered monomials ${\bf (Z37)}\,\,x^m\gL^n\,\,
(m,n\in{\bf Z})$ form a basis of the algebra and
$p=i\gl^{-1}x^{-1}(q^{1/2}\gL-q^{-1/2}\gL^{-1})$ where $\gl=q-q^{-1}\,\,(q\ne 1)$.
Note from \eqref{1a} there results
\bq\label{2a}
px=i\gl^{-1}(q^{-1/2}\gL-q^{1/2}\bar{\gL});\,\,xp=i\gl^{-1}(q^{1/2}
\gL-q^{-1/2}\bar{\gL})
\end{equation} 
Then the algebra can be represented as ${\bf (Z38)}\,\,\mf{H}_q=$ the associative
algebra
freely generated by $p,x,\gL,x^{-1},\gL^{-1}$ and their conjugates, modulo the ideal
generated by the relations \eqref{1a}, $\bar{x}=x,\,\,\bar{p}=p$, and
$\bar{\gL}=\gL^{-1}$.
\\[3mm]\indent
At the algebra level a field is an element of the subalgebra generated by $x$ and
$x^{-1}$ completed by formal series, i.e. ${\bf (Z39)}\,\,f(x)\in
[[x,x^{-1}]]\equiv\mf{A}_x$.  Then one has $pf(x)=g(x)p-iq^{1/2}h(x)\gL$ from the
algebra where $g$ and $h$ can be computed from \eqref{1a}.  Now define the derivative
${\bf (Z40)}\,\,\na:\,\mf{A}_x\to\mf{A}_x$ by $\na f(x)=h(x)$.  Since the monomials 
$x^m\,\,m\in{\bf Z}$ and 
\bq\label{3a}
\na x^m=[m]x^{m-1};\,\,[m]=\frac{q^m-q^{-m}}{q-q^{-1}}=q^{m-1}+q^{m-3}+\cdots+q^{-m+1}
\end{equation}
we see that $x^{-1}$ is not in the range of $\na$ and $ker(\na)=$ constants.  Similarly
one can define maps $L,\,L^{-1}:\,\mf{A}_x\to\mf{A}_x$ (onto) as follows.
Algebraically one can write ${\bf (Z41)}\,\,\gL f(x)=j(x)\gL$ and
$\gL^{-1}f(x)=k(x)\gL^{-1}$ so define
\bq\label{4a}
Lf(x)=j(x);\,\,L^{-1}f(x)=k(x)\Rightarrow Lx^m=q^{-m}x^m;\,\,L^{-1}x^m=q^mx^m
\end{equation}
Now the elements $x,\,x^{-1}\in\mf{A}_x$ also define maps $\mf{A}_x\to\mf{A}_x$ in the
obvious manner and thereby form an algebra
\bq\label{5a}
Lx=q^{-1}xL;\,\,L\na=q\na L;\,\,q^{1/2}x\na-q^{-1/2}\na x=-q^{-1/2}L
\end{equation}
homomorphic to the algebra \eqref{1a} with the identifications ${\bf (Z42)}\,\,L\sim
\gL,\,\,x\sim x,\,\,-iq^{1/2}\na\sim p$.  Without any bar operators being defined on
$L$ or $\na$ one verifies directly from the definitions of $L,\,L^{-1},$ and $\na$ that
\bq\label{6a}
\na=\gl^{-1}x^{-1}(L^{-1}-L);\,\,\na x^m=\frac{1}{\gl}(q^m-q^{-m})x^{m-1}=[m]
x^{m-1}
\end{equation}
which agrees with $\na x^m$ from \eqref{3a}.
\\[3mm]\indent
Next one determines a Leibnitz rule for $\na$ depending on the actions of L and
$L^{-1}$ on the product of fields.  Thus it is easily checked that ${\bf (Z43)}\,\,
L(x^nx^m)=(Lx^m)(Lx^n)$ and $L^{-1}(x^nx^m)=(L^{-1}x^m)(L^{-1}x^n)$ leading to
$L(fg)=(Lf)(Lg)$ and $L^{-1}(fg)=(L^{-1}f)(L^{-1}g)$.
For maps $x,\,x^{-1}$ one has evidently ${\bf (Z44)}\,\,xfg=(xf)g\equiv f(xg)$ and
$x^{-1}fg=f(x^{-1}g)$ (think of formal power series) and one obtains a Leibnitz rule
for
$\na$, namely, via the two calculations
\bq\label{7a}
\na
fg=\gl^{-1}x^{-1}(L^{-1}-L)fg=\gl^{-1}x^{-1}\left((L^{-1}f)(L^{-1}g)-(Lf)(Lg)\right);
\end{equation}
$$\na fg
=\left(\gl^{-1}x^{-1}(L^{-1}-L)f\right)(L^{-1}g)+\gl^{-1}(x^{-1}Lf)(L^{-1}g)+$$
$$+(Lf)\gl^{-1}x^{-1}(L^{-1}-L)g-\gl^{-1}(Lf)(x^{-1}L^{-1})g$$
Since $f$ and $g$ commute there results
\bq\label{8a}
\na fg=(\na f)(L^{-1}g)+(Lf)(\na g)=(\na f)(Lg)+(L^{-1}f)(\na g)
\end{equation}
Further for a Green's theorem one computes
$$\na(\na f)(L^{-1}g)=(\na^2f)g+(L^{-1}\na f)(\na L^{-1}g);\,\,\na(L^{-1}f)(\na g)=(\na
L^{-1}f)(L^{-1}\na g)+f(\na^2g);$$
\bq\label{9a}
(\na^2f)(g)-(f)(\na^2g)=\na\left((\na f)(L^{-1}g)-(L^{-1}f)(\na g)\right)
\end{equation}
\indent
An indefinite integral is defined as the inverse of $\na$.  First note (recall
$x^{-1}$ is not in the range of $\na$ and $\na c=0$)
\bq\label{10a}
\int^x x^n=\frac{1}{[n+1]}x^{n+1}+c
\end{equation}
Using \eqref{6a} one has now formally
$\na^{-1}=\gl (L^{-1}-L)^{-1}$
(not defined on $x^{-1}$) and a simple check shows agreement with \eqref{10a} on $x^n$
(with $c=0$) leading to 
\bq\label{12a}
\na^{-1}f(x)=\gl\sum_0^{\infty}L^{2\nu}Lxf(x)=-\gl\sum_0^{\infty}L^{-2\nu}L^{-1}xf(x)
\end{equation}
One uses the first or second series depending on which converges.  Again this checks
on $f=x^n$.  Thus by definitions ${\bf (Z45)}\,\,\int^x\na f=f+c$ and combining this
with \eqref{8a} gives integration by parts in the form
\bq\label{13a}
\int^x\na fg=fg+c=\int^x (\na f)(L^{-1}g)+\int^x(Lf)(\na g)\,\,\,OR
\end{equation}
$$\int^x\na fg=fg+c=\int^x(\na f)(L(g)+\int^x(L^{-1}f)(\na g)$$
\indent
At this point in \cite{wzz} a time variable for fields is introduced in the algebra
with $\bar{t}=t$ to produce an algebra $\mf{A}_{x,t}$ in which Schr\"odinger and 
Klein-Gordon equations are defined, along with other field equations, leading to gauge
theories in a purely algebraic context (with covariant derivatives, exterior forms,
differentials, connections, curvature, etc.).  We postpone this momentarily in order
to describe q-Fourier transforms, integration, and the corresponding natural Hilbert
spaces of $L^2$ type.
\\[3mm]\indent
One is looking for good representation spaces for $\mf{A}_x$ where $x$ and $p$ can be
diagonalized and with explicit formulas between $x$ and $p$ bases.  A natural choice
is e.g. the $Sin_q$ and $Cos_q$ functions (or other special q-functions) and Fourier
transformations (or corresponding eigenfunction transformations), from which one can
hopefully produce transmutation kernels in the spirit of Sections 2 and 3.  In this
respect  we will naturally be looking for the natural distribution like
generalizations as we go along in order to achieve suitable pairings of
q-eigenfunctions etc.  We sketch first some of the q-Fourier theory from \cite{wzz}  
(cf. also \cite{kw}).  Define
\bq\label{14a}
cos_q(x)=\sum_0^{\infty}(-1)^k\frac{x^{2k}}{[2k]!}\frac{q^{-k}}{\gl^{2k}};\,\,
sin_q(x)=\sum_0^{\infty}(-1)^k\frac{x^{2k+1}}{[2k+1]!}\frac{q^{k+1}}{\gl^{2k+1}}
\end{equation}
($\gl=q-q^{-1}$).  These functions satisfy (note $1-q^{-2k(2k+1)}=(\gl/q)q^{-2k}[2k+1]$)
\bq\label{15a}
\frac{1}{x}\left(sin_q(x)-sin_q(q^{-2}x)\right)=cos_q(x);\,\,
\frac{1}{x}\left(cos_q(x)-Cos_q(q^{-2}x)\right)=-q^{-2}sin_q(q^{-2}x)
\end{equation}
and \eqref{15a} in fact determines $sin_q(x)$ and $cos_q(x)$ up to an overall
normalization.  Note that \eqref{15a} corresponds to the usual derivative formulas for
sine and cosine.  Now $sin_q(x)$ and $cos_q(x)$ each form a complete and orthogonal
set in the following sense.  One defines a q-Fourier transform for suitable functions
$g(q^{2n})$ defined on lattice points $q^{2n}\,\,(n\in{\bf Z})$ via
\bq\label{16a}
\tl{g}_c(q^{2\nu})=N_q\sum_{-\infty}^{\infty}q^{2n}cos_q(q^{2(\nu+n)})g(q^{2n});
\end{equation}
$$g(q^{2n})=N_q\sum_{-\infty}^{\infty}q^{2\nu}cos_q(q^{2(\nu+n})\tl{g}_c(q^{2\nu})$$
with ($q>1$)
\bq\label{17a}
\sum_{-\infty}^{\infty}q^{2n}|g(q^{2n})|^2=\sum_{-\infty}^{\infty}
q^{2\nu}|\tl{g}_c(q^{2\nu})|^2;\,\,
N_q=\prod_0^{\infty}\left(\frac{1-q^{-2(2\nu+1)}}{1-q^{-4(\nu+1)}}\right)
\end{equation}
Similarly
\bq\label{18a}
\tl{g}_s(q^{2\nu})=N_q\sum_{-\infty}^{\infty}q^{2n}sin_q(q^{2(\nu+n)})g(q^{2n});
\end{equation}
$$g(q^{2n})=N_q\sum_{-\infty}^{\infty}q^{2\nu}
sin_q(q^{2(\nu+n)})\tl{g}_s(q^{2\nu});\,\,
\sum_{-\infty}^{\infty}q^{2n}|g(q^{2n})|^2=\sum_{-\infty}^{\infty}q^{2\nu}|\tl{g}_s
(q^{2\nu})|^2$$
(note $x=q^{2n}\geq 0$).  Further
\bq\label{19a}
N_q^2\sum_{-\infty}^{\infty}q^{2\nu}\left\{\begin{array}{c}
cos_q\\
sin_q\end{array}\right\}             
(q^{2(n+\nu)})\left\{\begin{array}{c}
cos_q\\
sin_q\end{array}\right\}
(q^{2(m+\nu)})=q^{-2n}\gd_{nm}
\end{equation}
and one has relations ${\bf (Z46)}\,\,cos_q(x)cos_q(qx)+q^{-1}sin_q(x)sin_q(q^{-1}x)=1$
(cf. \cite{wzz} for details).  From \eqref{19a} we see that $cos_q(q^{2n})$ and
$sin_q(q^{2n})$ must tend to zero for $n\to\infty$.  However one notes that 
$cos_q(q^{2n+1})$ and $sin_q(q^{2n+1})$ diverge for $n\to \infty$ and thus, although
$cos_q(x)$ and $sin_q(x)$ diverge for $x\to\infty$, the points $x=q^{2n}$ are close
to the zeros of $cos_q(x)$ and $sin_q(x)$ and for $n\to\infty$ tend to these zeros
such that the sum in \eqref{19a} is convergent.  One can also consider these functions
as a field, i.e. as elements of $\mf{A}_x$, and apply $\na$ in the form \eqref{6a} to 
get
\bq\label{20a}
\na cos_q(kx)=\frac{1}{\gl}\frac{1}{x}\{cos_q(qkx)-cos_q(q^{-1}kx)\}
\end{equation}
and setting $y=qkx$ in \eqref{15a}
this leads to ${\bf (Z47)}\,\,\na cos_q(kx)=-k(1/q\gl)sin_q
(q^{-1}kx)$; similarly $\na sin_q(kx)=k(q/\gl)cos_q(qkx)$.  This shows that
\bq\label{21a}
\na^2cos_q(kx)=-\frac{k^2}{q\gl^2}cos_q(kx);\,\,\na^2sin_q(kx)=
-\frac{k^2q}{\gl^2}sin_q(kx)
\end{equation}
and one notes also (cf. \cite{cqc})
\bq\label{211a}
\na cos_q(x)=-\frac{1}{q\gl}L sin_q(x);\,\,\na sin_q(x)=\frac{q}{\gL}L^{-1}cos_q(x)
\end{equation}
\indent
In order to produce an analogue of distribution pairings which were used extensively
in developing the operator theory of Section 3 one should now think of Fourier
transforms such as \eqref{16a} but defined on a specific class of functions $g$
corresponding to ${\mc D}\sim C_0^{\infty}$ or the Schwartz space ${\mc S}$ 
of rapidly decreasing functions and then
extend matters by duality as in distribution theory.  The basics here have already been
developed in \cite{oe} and are used below (cf. Section 5).
We omit here discusions of representations and $L^2$ spaces from \cite{wzz}
(which becomes quite complicated); this will be
examined in \cite{ca}.

\subsection{{\bf Heisenberg in higher dimensions.}}

Now for quantum groups and the R matrix we follow \cite{wzz} (cf. also
\cite{ca,cqc,cqd,ma,maa,se,sw,sh,wd}) and  
we will expound at more length in \cite{ca}.
As a 2-dimensional model consider matrices
\bq\label{44a}
T=\left(\begin{array}{cc}
a & b\\
c &d
\end{array}\right);\,\,ab=qba,\,\,ac=qca,\,\,ad=da+\gl
bc,\,\,bc=cb,\,\,bd=qdb,\,\,cd=qdc
\end{equation}
where $q\in {\bf C},\,\,q\ne 0,$ and $\gl=q-q^{-1}$.  Note $det_qT=ad-qbc$ is central
(i.e. commutes with $a,b,c,d$).  One considers the free associative algebra generated
by $1,a,b,c,d$ modulo the ideal of relations \eqref{44a} and in this algebra  
$\mf{A}$ formal power series are also allowed.  If $det_qT\ne 0$ one says 
$T\in GL_q(2)$ and if $det_qT=1$ then $T\in SL_q(2)$.  This is equivalent to other
authors with $q\leftrightarrow q^{-1}$.  The relations have very nontrivial
consequences.  In particular they allow an ordering of the elements $a,b,c,d$.  One
could decide to order monomials of degree n via a choice
$a^kb^{\ell}c^md^p$ for example with $n=k+{\ell}+m+p$.  Then it turns out that the
monomials of a given degree with such ordering are a basis for polynomials of fixed
degree (Poincar\'e- Birkhoff-Witt = PBW).  That the algebra $\mf{A}$ has the PBW
property follows from the fact that it can be formulated with the help of an R matrix. 
Thus the relations \eqref{44a} can be written in the form ${\bf (Z48)}\,\,\sum
R^{ij}_{k\ell} T^k_rT^{\ell}_s=\sum T^i_kT^j_{\ell}R^{k\ell}_{rs}$ where the indices
take values 1 and 2.  R can be written in the form
\bq\label{45a}
R=\left(\begin{array}{cccc}
q & 0 & 0 & 0\\
0 & \gl & 1 & 0\\
0 & 1 & 0 & 0\\
0 & 0 & 0 &  q
\end{array}\right)
\end{equation}
The rows and columns are labelled by $11,\,12,\,21,\,$ and $22$.  As an example
consider
\bq\label{46a}
\sum R^{12}_{ij}T^i_2T^j_2=\sum T^1_iT^2_jR^{ij}_{22}\equiv \gl T^1_2T^2_2+T^2_2T^1_2=
qT^1_2T^2_2
\end{equation}
which is $\gl bd +db=qbd\Rightarrow bd=qdb$.  The relations {\bf (Z48)} are called
RTT relations; there are 16 of these reducing to the 6 relations of \eqref{44a}.  Since
R matrices are a defining property of quasitriangular Hopf algebras or quantum groups
leading to braiding etc. the PBW property seems to be another way to characterize this.
In any case there are far reaching consequences; e.g. from the RTT relations follows
(cf. \cite{wzz} for details) ${\bf (Z49)}\,\,\gD T^j_{\ell}=\sum T^j_r\ot T^r_{\ell}$
which is compatible with the RTT relations
${\bf (Z50)}\,\,\sum R^{ij}_{k\ell}\gD T^k_r\gD T^{\ell}_s=\sum \gD
T^i_k\gD T^j_{\ell}R^{k\ell}_{rs}$.  Further for an antipode S one enlarges the algebra
by the inverse of $det_qT$ to obtain ($S(T)\sim T^{-1}$ - cf. \cite{ca,sh})
\bq\label{48a}
T^{-1}=\frac{1}{det_qT}\left(\begin{array}{cc}
d & -q^{-1}b\\
-qc & a
\end{array}\right)
\end{equation}
For the counit one takes $\gep(T)=id$ and we have a quasitriangular Hopf algebra.
For more details in this particular instance note the characteristic equation
${\bf (Z51)}\,\,(R-q)(R+q^{-1})=0$ where $q,\,-q^{-1}$ are the eigenvalues for R.
The projectors onto the respective eigenspaces are
\bq\label{49a}
A=-\frac{1}{1+q^2}(R-q);\,\,B=\frac{q}{1+q^2}(R+q^{-1})
\end{equation}
A is a deformation of an antisymmetrizer and B of a symmetrizer.  The normalization is
such that 
\bq\label{50a}
A^2=A,\,\,B^2=B,\,\,AB=BA=0,\,\,1=A+B,\,\,R=qB-q^{-1}A
\end{equation}
\begin{example}
This approach can be generalized to n dimensions with an $n^2\times n^2$ matrix for
$GL_q(N)$ of the form
\bq\label{51a}
R^{ji}_{k\ell}=\gd_k^i\gd_{\ell}^j
[1+(q-1)\gd^{ij}]+(q-q^{-1})\gt(i-j)\gd_k^j\gd^i_{\ell}
\end{equation}
where $\gt(i-j)=1$ for $i>j$ and zero otherwise.  The RTT relations {\bf (Z48)} now
refer to an $n\times n$ matrix T and there are $n^4$ relations for the $n^2$ entries
of T.  It can be shown that the T matrix so defined has the PBW property. 
Comultiplication is defined as in {\bf (Z49)} and the R matrix satisfies the same
characteristic equation {\bf (Z51)}.  These R matrices \eqref{45a} and \eqref{51a} are
symmetric ${\bf (Z52)}\,\,R^{ab}_{cd}=R^{cd}_{ab}$ and in such a situation the
transposed matrix $\tl{T}:\,\,\tl{T}^a_b=T^b_a$ also satisfies the RTT relations
{\bf (Z48)} (cf. \cite{wzz} for a calculation).  Further for $T\in GL_q(n)$ one has
$\tl{T}\in GL_q(n)$ but (for $q\ne 1$) $\tl{T}^{-1}\ne \widetilde{T^{-1}}$.
\end{example}
\indent
Given now some quantum groups one looks at their comodules which are referred to as
quantum planes (cf. also \cite{ma}).  We give results here only for two dimensional
situations and refer to \cite{ca,wzz} for the general situation.
For $GL_q(n)$ with $A\sim R-q$ (cf. \eqref{49a}) one has ${\bf (Z56)}\,\,
x^ix^j=(1/q)\sum R^{ij}_{k\ell}x^kx^{\ell}$ and in two dimensions this reduces to the
condition ${\bf (Z57)}\,\,x^1x^2=qx^2x^1$.  The relation {\bf (Z56)} can be
generalized to the situation of two or more copies of quantum planes, e.g. $(x^1,x^2)$
and $(y^1,y^2)$ and the relations ${\bf (Z58)}\,\,x^iy^j=(\kappa/q)\sum
R^{ij}_{k\ell}y^kx^{\ell}$ are consistent (i.e. they have the PBW property for
arbitrary $\kappa\ne 0$ and they are covariant).  For $n=2$ {\bf (Z58)} becomes
\bq\label{53a}
x^1y^1=\gk y^1x^1,\,\,x^1y^2=\frac{\gk}{q}y^2x^1+\frac{\gk\gl}{q} y^1x^2,\,\,
x^2y^1=\frac{\gk}{q}y^1x^2,\,\,x^2y^2=\gk y^2x^1
\end{equation}
Consistency can be checked directly (cf. \cite{wzz}).  One can also conclude that 
\eqref{53a} does not generate new relations and therefore the PBW property holds.
This can be done more systematically by starting from a general R matrix and
considering three copies of quantum planes, $x,\,y,\,z$; the covariant relations are
${\bf (Z59)}\,\,xy=Ryx,\,\,yz=Rzy,\,\,xz=Rzx$ (indices as in {\bf (Z58)}).  One
demands that a reordering of $xyz$ to $zyx$ should give the same result, independently
of the way in which this reordering is achieved.  There are two independent ways to do
this and their equivalence leads to the quantum Yang Baxter equation (QYBE) which can
be formulated by introducing $n^3\times n^3$ matrices
${\bf (Z60)}\,\,R_{12j}^i=R^{i_1i_2}_{j_1j_2}\gd^{i_3}_{j_3}$ and $R_{23j}^i=
\gd^{i_1}_{j_1}R^{i_2i_3}_{j_2j_3}$ and the QYBE is then ${\bf
(Z61)}\,\,R_{12}R_{23}R_{12}=R_{23}R_{12}R_{23}$.  There are $n^6$ equations for $n^4$
independent entries of the R matrix.  It can be checked that the matrices \eqref{45a}
and \eqref{51a} satisfy the QYBE. 
\\[3mm]\indent
Another algebraic structure on comodules is obtained by generalizing the Leibnitz rule
${\bf (Z65)}\,\,(\pp/\pp x^i)x^j=\gd_i^j+x^j(\pp/\pp x^i)$.  One demands that the
algebra generated by the elements of the quantum plane algebra $x^i$ and the
derivatives $\pp_i$, modulo proper ideals, have the PBW property.  In addition one
shows that there is an exterior DC based on these quantum properties.  Beginning with an 
Ansatz ${\bf (Z66)}\,\,\pp_ix^j=\gd_i^j+\sum C^{jk}_{i\ell}x_k\pp^{\ell}$.  
and using QYBE and covariance ideas one arrives at derivative rules 
consistent with the quantum plane relations.  For
$SL_q(2)$ one has ${\bf (Z75)}\,\,R=qB-q^{-1}A$ and
$R^{-1}=q^{-1}B-qA$ (cf. \eqref{50a})
which implies $C=qR$ or $C^{-1}=q^{-1}R^{-1}$.  There are then
two solutions with the desired properties and one concludes that
\bq\label{63a}
\pp_ix^j=\gd_i^j+\sum
qR^{jk}_{i\ell}x^{\ell}\pp_k\,\,{\bf or}\,\,\,\hat{\pp}_ix^j=\gd_i^j+q^{-1}\sum
(R^{-1})^{jk}_{i\ell}x^{\ell}\hat{\pp}_k
\end{equation} 
are two possibilities to define a covariant derivative on a quantum plane consistent
with the defining relations of $SL_q(n)$ quantum planes.  
In summary one has for the algebra based on the R matrix 
\eqref{45a} for $SL_q(2)$
\begin{enumerate}
\item[3.]
$x^ix^j=q^{-1}\sum R^{ij}_{k\ell}x^kx^{\ell};\,\,x^1x^2=qx^2x^1$
\item[4.]
$\pp_ix^j=\gd_i^j+qR_{i\ell}^{jk}x^{\ell}\pp_k;\,\,\pp_1x^1=1+q^2x^1\pp_1+q\gl
x^2\pp_2;\,\,\pp_1x^2=qx^2\pp_1;\,\,\pp_2x^1=qx^1\pp_2;\,\,\pp_2x^2=1+q^2x^2\pp_2$
\item[5.]
$\pp_a\pp_b=q^{-1}\sum\pp_c\pp_dR^{dc}_{ba};\,\,\pp_1\pp_2=q^{-1}\pp_2\pp_1$
\item[6.]
$\hat{\pp}_ix^j=\gd_i^j+q^{-1}\sum
(R^{-1})^{jk}_{i\ell}x^k\hat{\pp}_{\ell};\,\,\hat{\pp}_1x^1=1+q^{-2}x^1\hat{\pp}_1;\,\,
\hat{\pp}_1x^2=q^{-2}x^2\hat{\pp}_1;\,\,
\hat{\pp}_2x^1=q^{-1}x^1\hat{\pp}_2;\\1+q^{-2}x^2\hat{\pp}_2-
\gl q^{-1}x^1\hat{\pp}_1=\hat{\pp}_2x^2$
\item[7.]
$\hat{\pp}_a\hat{\pp}_b=q^{-1}\sum\hat{\pp}_c\hat{\pp}_dR^{dc}_{ba};\,\,
\hat{\pp}_1\hat{\pp}_2=q^{-1}\hat{\pp}_2\hat{\pp}_1;\,\,\hat{\pp}_a\pp_b=
q\sum R^{cd}_{ba}\pp_d\hat{\pp}_c;\,\,\hat{\pp}_1\pp_1=q^2\pp_1\hat{\pp}_1;\,\,
\hat{\pp}_1\pp_2=q\pp_2\hat{\pp}_1;\,\,\hat{\pp}_2\pp_1=q\pp_1\hat{\pp}_2+\gl
q\pp_2\hat{\pp}_1;\,\,\hat{\pp}_2\pp_2=q^2\pp_2\hat{\pp}_2$
\end{enumerate}
\indent
There is also an exterior DC based on these quantum derivatives involving
\bq\label{66a}
d^2=0,\,\,d(fg)=(df)g+fdg,\,\,ddx^i=-dx^id
\end{equation}
\indent
We can now deal with the entire algebra generated by $x^i,\,dx^j,$ and $\pp_{\ell}$
modulo the respective ideals.  For this purpose the $dx^{\ell},\,\pp_j$ relations have
to be specified and we provide he $xdx$ relations for $SL_q(2)$ via 
\bq\label{71a}
dx^1x^1=q^2x^1dx^1;\,\,dx^1x^2=qx^2dx^1+(q^2-1)x^1dx^2;
\end{equation}
$$dx^2x^1=qx^1dx^2;\,\,dx^2x^2=q^2x^2dx^2$$
Next one adds a conjugation via ${\bf (W1)}\,\,\ol{x^i}=\bar{x}^i$ and
$\ol{x^ix^j}=\bar{x}_j\bar{x}_i$ (the lower index is for convenience) and
for $n=2$ the explicit $x\bar{y}$ relations are (take $\gk=1$)
\bq\label{74a}
x^1\bar{y}_1=\bar{y}_1x^1-q\gl\bar{y}_2x^2;\,\,x^1\bar{y}_2=q\bar{y}_2x^1;\,\,x^2\bar{y}_1=
q\bar{y}_1x^2;\,\,x^2\bar{y}_2=\bar{y}_2x^2
\end{equation}
The $\bar{y}\bar{x}$ relations follow from \eqref{53a} by conjugation.  For the
entries of T and $\bar{T}$ as defined by \eqref{44a} and its conjugate the relations
are
\bq\label{75a}
a\bar{a}=\bar{a}a-q\gl\bar{c}c;\,\,a\bar{b}=q^{-1}\bar{b}a-\gl\bar{d}c;\,\,a\bar{c}
=q\bar{c}a;\,\,a\bar{d}=\bar{d}a;
\end{equation}
$$b\bar{a}=q^{-1}\bar{a}b-\gl\bar{c}d;\,\,b\bar{b}=\bar{b}b+q\gl(\bar{a}a-
\bar{d}d-q\gl\bar{c}c);\,\,b\bar{c}=\bar{c}b;$$
$$b\bar{d}=q\bar{d}b=\gl
q^2\bar{c}a;\,\,c\bar{a}=q\bar{a}c;\,\,c\bar{b}=\bar{b}c;\,\,c\bar{c}=\bar{c}c;\,\,
c\bar{d}=q^{-1}\bar{d}c;$$
$$d\bar{a}=\bar{a}d;\,\,d\bar{b}=q\bar{b}d+\gl q^2\bar{c}a;\,\,d\bar{c}=q^{-1}\bar{c}d;
\,\,d\bar{d}=\bar{d}d+\gl q\bar{c}c$$
It is possible to identify
$\bar{T}$ with $T^{-1}$ and we get the quantum group $U_q(n)$ or for $det_qT=1$ the
quantum group $SU_q(n)$.
For $n=2$ again
\bq\label{76a}
T=\left(\begin{array}{cc}
a & b\\
c & d
\end{array}\right);\,\,T^{-1}=\left(\begin{array}{cc}
d & -q^{-1}b\\
-qc & a
\end{array}\right)
\end{equation}
and one finds ${\bf (W9)}\,\,\bar{a}=d,\,\,\bar{b}=-qc,\,\,\bar{d}=a,$ and
$\bar{c}=-q^{-1}b$.  It can be verified directly that {\bf (W9)} is consistent with 
\eqref{75a}.  

\subsection{{\bf Q-deformed Heisenberg algebras.}}

The canonical commutation relations here are ${\bf (W12)}\,\,[\hat{x},\hat{p}]=i$ with
self adjoint $\hat{x}$ and $\hat{p}$ and a physical system is defined via a
representation of this algebra in a Hilbert space.  In $L^2$ like spaces one has
$\hat{x}=x$ and $\hat{p}\sim -i(\pp/\pp x)$.  Now {\bf (W12)} will be changed in accord
with quantum group considerations.  It is natural to assume that $\hat{x}$ is an
element of a quantum plane and $\hat{p}$ should be a derivative in this plane; the
simplest example is suggested by equations in (3)-(7) in Section 4 above and one is
led to consider the algebra involving ${\bf (W13)}\,\,\pp x=1+qx\pp$.  More precisely
one looks at the algebra generated by the elements x and $\pp$ modulo the relation
{\bf (W13)}.  If one assumes $\bar{x}=x$ then $i\pp$ cannot be selfadjoint since
from {\bf (W13)} one has ${\bf (W14)}\,\,\bar{\pp}x=-q^{-1}+q^{-1}x\bar{\pp}$ and in
general $\bar{\pp}$ will be related to $\hat{\pp}$ rather than to $\pp$ (cf. (3)-(7)
again).  Thus one could study the algebra generated by $x,\,\pp,$ and $\bar{\pp}$
modulo {\bf (W13)}, {\bf (W14)}, and an ideal generated by $\pp\bar{\pp}$ relations. 
These can be found by an argument as before (cf. {\bf (Z81)}) and using {\bf (W13)},
{\bf (W14)} one arrives at ${\bf (W15)}\,\,\bar{\pp}=q\pp\bar{\pp}$.  This is
consistent but if we try to define an operator $\hat{p}=-(i/2)(\pp-\bar{\pp})$ the 
$x,\,\hat{p}$ relations do not close.  However $\bar{\pp}$ can be related to $\pp$ and
x in a nonlinear way via the scaling operator $\gL$ and one has
\bq\label{79a}
\gL=q^{1/2}(1+(q-1)x\pp);\,\,\gL x=qx\gL;\,\,\gL\pp=q^{-1}\pp\gL
\end{equation}
The scaling property follows from {\bf (W13)}.  We now define
${\bf(W16)}\,\,\tl{\pp}=-q^{-1/2}\gL^{-1}\pp$ and $\gL^{-1}$ is defined by an
expansion in $(q-1)$ leading to 
\bq\label{80a}
\tl{\pp}x=-q^{-1}+q^{-1}x\tl{\pp};\,\,\tl{\pp}\pp=q\pp\tl{\pp}
\end{equation}
Comparing this with {\bf (W14)} and {\bf (W15)} it follows from {\bf (W16)} and 
\eqref{80a} that conjugation in the $x,\,\pp$ algebra can be defined by
${\bf (W17)}\,\,\bar{x}=x,\,\,\bar{\pp}=-q^{-1/2}\gL^{-1}\pp\sim\tl{\pp}$. 
Conjugating $\gL$ and using {\bf (W17)} shows that ${\bf (W18)}\,\,\bar{\gL}=\gL^{-1}$
so $\gL$ is unitary, which justifies the factor $q^{1/2}$ in the definition of $\gL$.
The existence of a scaling operator $\gL$ and the definition of conjugation via ${\bf
(W17)}$ seems to be very specific for the $x,\,\pp$ algebra {\bf (W13)}.  It is
however generic in the sense that a scaling operator and a definition of conjugation
based on it can be found for all the quantum planes defined by $SO_q(n)$ and
$SO_q(1,n)$.  The definition of the q-deformed Heisenberg algebra will now be based on
the definition of the momentum given by $p=-(i/2)(\pp-\bar{\pp})$ which is selfadjoint.
There results then
\bq\label{81a}
q^{1/2}xp-q^{-1/2}px=i\gL^{-1};\,\,\gL x=qx\gL;\,\,\gL p=q^{-1}p\gL;\,\,\bar{p}=p;\,\,
\bar{x}=x;\,\,\bar{\gL}=\gL^{-1}
\end{equation}
These algebraic relations can be verified in the $x,\,\pp$ representation where the
ordered $x,\,\pp$ monomials form a basis.  One takes \eqref{81a} as the defining
relations for the q-deformed Heisenberg algebra without making further reference to
its $x,\,\pp$ representation.
\\[3mm]\indent
A nice continuation of the Wess-Zumino calculus \cite{wd} is developed in \cite{oh}. 
Assume the quadratic commutation relations for the coordinates are given by ${\bf
(W19)}\,\,\sum P^{ij}_{k\ell}x^kx^{\ell}=0$ where $R=\ga P+R'$ with P a projection
(cf. \cite{rz,wd}).  Then the action of the derivatives is given by
\bq\label{4z}
\pp_ix^j=\gd^j_i-\ga^{-1}\sum R^{jk}_{i\ell}x^{\ell}\pp_k
\end{equation}
One checks that this formula is covariant under the corresponding quantum group 
(e.g. $GL_q(2)$) and consistent with the relations {\bf (Z3)}.  The isomorphisms to
follow now are based on the decomposition of the ring of q-differential operators into
the tensor product of rings $Diff_{q^2}(1)$ where $Diff_{q^2}(1)$ is the ring
generated by $x$ and $\pp^q$ obeying the relation ${\bf (W20)}\,\,\pp^qx=1+q^2x\pp^q$.
An isomorphism of $Diff_{q^2}(1)$ with the ring of usual differential operators in
one variable can be attained via
\bq\label{5z}
\pp^q=\frac{1}{x}\frac{e^{2hx\pp}-1}{q^2-1};\,\,q=e^h;\,\,\pp=\frac{\pp}{\pp x}
\end{equation}
Note here $exp(2hx\pp)f(x)=f(x+2hx)=f((1+2h)x)\sim f(q^2x)$ since $q^2=exp(2h)\sim 1+2h$
but \eqref{5z} is not the same as ${\bf (W21)}\,\,D_{q^2}f(x)=(f(q^2x)-f(x))/(x(q^2-1))$
and {\bf (W21)} will seem to be a more natural isomorphism later.
In the GL case the commutation relations for the coordinates are ${\bf (W22)}\,\,
x^ix^j=qx^jx^i\,\,(i<j)$ and the scheme above gives then
\bq\label{6z}
\pp_ix^j=qx^j\pp_i\,\,(i\ne j);\,\,\pp_ix^i=1+q^2x^i\pp_i+q\gl\sum_{j>i}x^j\pp_j
\end{equation}
where $\gl=q-q^{-1}$.  Now introduce the quantities ${\bf
(W23)}\,\,\mu_k=1+q\gl\sum_{j\geq k}x^j\pp_j$ and then one can rewrite \eqref{6z} in
the form ${\bf (W24)}\,\,\pp_ix^i=\mu_i+x^i\pp_i$.  If $\mu_i$ were a number one could
remove it by renormalizing the derivative but since it is an operator one can use the
following simple commutation relations
\bq\label{7z}
\mu_ix^j=x^j\mu_i\,\,(i>j);\,\,\mu_ix^j=q^2x^j\mu_i\,\,(i\leq j);
\end{equation}
$$\mu_i\pp_j=
\pp_j\mu_i\,\,(i>j);\,\,\mu_i\pp_j=q^{-2}\pp_j\mu_i\,\,(i\leq j)$$
This implies that ${\bf (W25)}\,\,\mu_i\mu_j=\mu_j\mu_i$.  Note that in the one
dimensional example one could write $\mu=1+q\gl x\pp^q$ for the operators satisfying 
{\bf (W20)} and then $\mu x=q^2x\mu$ with $\pp^q\mu=q^2\mu\pp^q$.  In terms of classical
variables one would have $\mu=exp(2hx\pp)$ which explains the multiplicative nature of
$\mu$.  Now one can define ${\bf (W26)}\,\,X^i=(\mu_i)^{-1/2}x^i$ and
$D_i=q(\mu_i)^{-1/2}\pp_i$.  Using the multiplicative properties \eqref{7z} one can
rewrite {\bf (W24)} in the new variables {\bf (W26)} now in a form which does not
contain $\mu_i$.  Actually the transformation {\bf (W26)} gives even more, namely
\bq\label{8z}
X^iX^j=X^jX^i;\,\,D_iD_j=D_jD_i;\,\,D_iX^j=X^jD_i\,\,(i\ne j);\,\,D_jX^j=1+q^{-2}X^jD_j
\end{equation}
Thus the relations \eqref{6z} are completely untangled by the transformation {\bf
(W26)} and the whole ring of q-differential operators becomes just the tensor product
of rings $Diff_{q^{-2}}(1)$ leading to the result that $Diff_{GL_q}(n)\simeq Diff(n)$
which is the ring of classical differential operators in n variables.
Note also $Diff_{q^{-2}}(1)\simeq Diff_{q^2}(1)$ via $\gd^q=q\mu^{-1/2}\pp^q$ and
$y=\mu^{-1/2}x$ leading to $\gd^qy=1+q^{-2}y\gd^q$.

\section{FOURIER TRANSFORMS AND DISTRIBUTIONS}
\renewcommand{\theequation}{5.\arabic{equation}}
\setcounter{equation}{0}

We go here to \cite{oe,rc} (cf. also \cite{kb,ki,kw,sz}) where a Fourier theory is
developed along the lines of ``classical" distribution theory (cf. \cite{cv,czc}).  One
deals here with a $q^2$ theory which in view of constructions in \cite{wzz} (cf. Section
4) should be easier to control.  A rather full exposition is given here (from \cite{oe})
in view of the importance of this subject.  Thus define (cf. \eqref{41b}-\eqref{43b})
\bq\label{1e}
\nm_{q^2}=\frac{(q^2;q^2)_n}{(q^2;q^2)_m(q^2;q^2)_{n-m}};\,\,(a;q)_n=\left\{
\begin{array}{cc}
1 & n=0\\
\prod_0^{n-1}(1-aq^j) & n\geq 1
\end{array}\right.
\end{equation}
and for $|z|<1$
\bq\label{2e}
e_{q^2}(z)=\sum_0^{\infty}\frac{z^n}{(q^2;q^2)_n}=\frac{1}{(z;q^2)_{\infty}};\,\,
E_{q^2}(z)=\sum_0^{\infty}\frac{q^{n(n-1)}z^n}{(q^2;q^2)_n}=(-z;q^2)_{\infty}
\end{equation}
\bq\label{3e}
cos_{q^2}z=\frac{1}{2}[e_{q^2}(iz)+e_{q^2}(-iz)];\,\,sin_{q^2}z=\frac{1}{2i}
[e_{q^2}(iz)-e_{q^2}(-iz)];
\end{equation}
$$Cos_{q^2}z=\frac{1}{2}[E_{q^2}(iz)+E_{q^2}(-iz)];\,\,Sin_{q^2}z=\frac{1}{2i}
[E_{q^2}(iz)-E_{q^2}(-iz)]$$
There is also a basic hypergeometric series ${\bf
(W27)}\,\,{}_0\Phi_1(-;0;q^2;z)=\sum_0^{\infty}[q^{2n(n-1)}z^n/(q^2;q^2)_n]$ and one shows
that
\bq\label{4e}
e_{q^2}(z)=\frac{1}{(q^2;q^2)_{\infty}}\sum_0^{\infty}
\frac{(-1)^kq^{k(k+1)}}{(q^2;q^2)_k(1-zq^{2k})}
\end{equation}
leading to
\bq\label{5e}
cos_{q^2}z=\frac{1}{(q^2;q^2)_{\infty}}
\sum_0^{\infty}\frac{(-1)^kq^{k(k+1)}}{(q^2;q^2)_k(1+z^2q^{4k})};
\end{equation}
$$sin_{q^2}z=\frac{z}{(q^2;q^2)_{\infty}}\sum_0^{\infty}
\frac{(-1)^kq^{k(k+3)}}{(q^2;q^2)_k(1+z^2q^{4k})}$$
It follows that for $z,q$ real
\bq\label{6e}
|cos_{q^2}z|\leq\frac{(-q^2;q^2)_{\infty}}{(1+z^2)(q^2;q^2)_{\infty}};\,\,
|sin_{q^2}z|\leq \frac{|z|(-q^2;q^2)_{\infty}}{(1+z^2)(q^2;q^2)_{\infty}}
\end{equation}
Further ${\bf (W28)}\,\,|Cos_{q^2}z|\leq 1$ and $|Sin_{q^2}z|\leq |z|$ and one has use for
the functions
\bq\label{7e}
{\bf Q}(z,q)=(1-q^2)\sum_{-\infty}^{\infty}\frac{1}{zq^{2m}+z^{-1}q^{-2m}};\,\,
\gT(z)=(1-q^2)\sum_{-\infty}^{\infty}sin_{q^2}((1-q^2)q^{2m}z)
\end{equation}
satisfying $(\bs)\,\,\gT(q^{2k}z)=\gT(z)$ for $z\ne 0$ and $\gT(z)={\bf Q}((1-q^2)z,q)$.
In this respect one has
\bq\label{8e}
\gT(z)=\frac{1-q^2}{(q^2;q^2)_{\infty}}\sum_{k=0}^{\infty}
\frac{(-1)^kq^{k(k+3)}}{(q^2;q^2)_k}\sum_{m=-\infty}^{\infty}
\frac{(1-q^2)zq^{2m}}{1+(1-q^2)^2z^2q^{4(m+k)}}=
\end{equation}
$$=\frac{1-q^2}{(q^2;q^2)_{\infty}}\sum_{k=0}^{\infty}
\frac{(-1)^kq^{k(k+1)}}{(q^2;q^2)_k}\sum_{-\infty}^{\infty}
\frac{(1-q^2)zq^{2m}}{1+(1-q^2)^2z^2q^{4m}}={\bf Q}((1-q^2)z,q)$$
(note also ${\bf (W29)}\,\,\gT_0=\gT(1)={\bf Q}(1-q^2,q)$).  Then for an arbitrary integer
$M>0$ one has
\bq\label{9e}
(1-q^2)z\sum_{m=-M}^{\infty}q^{2m}cos_{q^2}((1-q^2)q^{2m}z)=sin_{q^2}((1-q^2)q^{-2M}z);
\end{equation}
$$(1-q^2)z\sum_{m=-M}^{\infty}
q^{2m}sin_{q^2}((1-q^2)q^{2m}z)=1-cos_{q^2}((1-q^2)q^{-2M}z)$$
One defines now for $f\in\mf{A}={\bf C}(z,z^{-1})$ (formal Laurent series) 
$\pp_zf(z)=\{z^{-1}/(1-q^2)\}[f(z)-f(q^2z)]$ leading to
\bq\label{10e}
\pp_z^kz^n=\left\{\begin{array}{cc}
\frac{(q^2;q^2)_n}{(q^2;q^2)_{n-k}(1-q^2)^k}z^{n-k} & 0\leq k\leq n\\
0 & k>n
\end{array}\right.
\end{equation}
$$\pp_z^kz^{-n-1}=(-1)^kq^{-k(2n+k+1)}\frac{(q^2;q^2)_{n+k}}{(q^2;q^2)_n(1-q^2)^k}
z^{-n-k-1}$$
Now the $q^2$-Jackson integral is defined as (cf. Section 6 for some background)
\bq\label{11e}
I_{q^2}f=\int d_{q^2}zf(z)=(1-q^2)\sum_{-\infty}^{\infty}q^{2m}[f(q^{2m})+f(-q^{2m})]
\end{equation}
\begin{definition}
The function $f(z)$ is locally $q^2$ integrable if the $q^2$ integral
\bq\label{12e}
\int_a^bd_{q^2}zf(z)=(1-q^2)\sum_0^{\infty}q^{2m}[bf(bq^{2m})-af(aq^{2m})]
\end{equation}
exists for any finite $a,b$ (i.e. the series in \eqref{12e} converges).  $f(z)$ is
absolutely $q^2$ interable if the series $\sum_{-\infty}^{\infty}q^{2m}[|f(q^{2m}|+
|f(-q^{2m})|]$ converges.
\end{definition}
\indent
Let now $\mf{B}$ be the algebra analogous to $\mf{A}$ but generated by $s,s^{-1}$ with
$zs=q^2sz$.  Define $q^2$ differentiation in $\mf{B}$ via ${\bf
(W30)}\,\,\pp_s\phi(s)=(\phi(s)-\phi(q^2s)]s^{-1}/(1-q^2)$.  One denotes by $\mf{A}\mf{B}$
the algebra with generators $z,z^{-1},s,s^{-1}$ and the $q^2$ differentiation
${\bf (W31)}\,\,zs=q^2sz,\,\,\pp_zs=q^{-2}s\pp_z,\,\,\pp_sz=q^2z\pp_s,$ and $\pp_z\pp_s=
q^2\pp_s\pp_z$.  One considers $\mf{A}\mf{B}$ as a left module under the left action of
$\mf{A}$ by multiplication and a right module under the action of $\mf{B}$.  In
order to define a $q^2$ integral on $\mf{A}\mf{B}$ one orders the generators of the
integrand so that z stays on the left while s stays on the right.  For example if $f(z)=
\sum a_rz^r$ then ${\bf (W32)}\,\,f(zs)=\sum a_r(zs)^r=\sum a_r q^{-r(r-1)}z^rs^r$.  For
convenience one writes ${\bf (W33)}\,\,\ddg g(zs)\ddg=\sum a_rz^rs^r$ if $g(z)=\sum
a_rz^r$.  Then for example one will have ${\bf (W34)}\,\,E_{q^2}(i(1-q^2)zs)=\ddg 
e_{q^2}(i(1-q^2)zs)\ddg$.  One can now calculate
\begin{itemize}
\item
${\bf (W35)}\,\,\int d_{q^2}zz^{-1}E_{q^2}(i(1-q^2)zs)=2i\gT (s)$ (see {\bf (W29)}
\item
${\bf (W36)}\,\,\int d_{q^2}z\gT(z)E_{q^2}(i(1-q^2)zs)=2i\gT_0s^{-1}$ (use \eqref{9e} and
\eqref{6e})
\end{itemize}
\indent
Let ${\mc S}_{q^2}=\{\phi(x)\}$ be the space of infinitely $q^2$ differentiable rapidly
decreasing functions ${\bf (W37)}\,\,|x^k\pp_x^{\ell}\phi(x)|\leq C_{k,\ell}(q)$ for
$k\geq 0,\,\ell\geq 0$.  Let ${\mc S}$ be the space of $C^{\infty}$ functions in the
classical sense (so ${\bf (W38)}\,\,|x^k\phi^{(\ell)}(x)|\leq C_{k,\ell}$ for $k,\ell\geq
0$.  Then ${\mc S}\subset{\mc S}_{q^2}$.  To see this let $\phi(x_1,x_2,\cdots,x_k)$ be
the separated difference of order k, so for an arbitrary integer $\ell\geq 0$
\bq\label{122e}
\phi(q^{2\ell}x,q^{2\ell-2}x,\cdots,x)=\frac{(q^2;q^2)_{\ell}}{(1-q^2)^{\ell}}\pp_x^{\ell}\phi(x)
\end{equation}
On the other hand if $\phi(x)$ is $\ell$ times differentiable in the classical sense one
has
\bq\label{13e}
\phi(q^{2\ell}x,q^{2\ell-2}x,\cdots,x)=\frac{1}{\ell!}\phi^{(\ell)}(\xi)\,\,\,(\xi\in
(q^{2\ell}x,x))
\end{equation}
It follows that
\bq\label{14e}
\pp_x^{\ell}\phi(x)=\frac{(1-q^2)^{\ell}}{(q^2;q^2)_{\ell}\ell!}\phi^{(\ell)}(\xi)
\end{equation}
and therefore from {\bf (W38)} 
\bq\label{15e}
|x^k\pp_x^{\ell}\phi(x)|\leq\frac{(1-q^2)^{\ell}}{(q^2;q^2)_{\ell}\ell!}C_{k,\ell}
\end{equation}
\begin{definition}
The skeleton $\hat{\phi}(z)$ of $\phi(z)\in {\mc S}_{^2}$ is the set of evaluations of 
$\phi(z)$ on the lattice generated by the powers of $q$, i.e. ${\bf (W39)}\,\,\hat{\phi}
(z)=\{\phi(q^{2n}),\,\,n=0,\pm 1,\cdots\}$.  The space of skeletons is denoted by
$\hat{{\mc S}}_{q^2}$ and a basis in this space is generated by
\bq\label{16e}
\hat{\phi}_n^{+}(z)=\left\{\begin{array}{cc}
1 & z=q^{2n}\\
0 & z\ne q^{2n}\end{array}\right.;\,\,\hat{\phi}_n^{-}(z)=\left\{\begin{array}{cc}
1 & z=-q^{2n}\\
0 & z\ne -q^{2n}
\end{array}\right.
\end{equation}
so that
$\hat{\phi}(z)=\sum_{-\infty}^{\infty}[a_n\hat{\phi}_n^{+}(z)+b_n\hat{\phi}_n^{-}(z)]$
for any $\hat{\phi}\in \hat{{\mc S}}_{q^2}$.  The topology in $\hat{{\mc S}}_{q^2}$ is
induced by that in ${\mc S}_{q^2}$
\end{definition}
\indent
Let $L_{q^2}$ be the linear map from ${\mc S}_{q^2}$ to $\hat{{\mc S}}_{q^2}$ defined by
the evaluation of functions at the vertices of the lattice.  Then ${\bf (W40)}\,\,
\hat{{\mc S}}_{q^2}={\mc S}_{q^2}/ker(L_{q^2})$.  Let $\gL\phi(z)=\phi(q^2z)$ and then
${\bf (W41)}\,\,\gL z=q^2z\gL$ with $\pp_z\gL=q^2\gL \pp_z$.  The operations $\gL$ and
$\pp_z$ are well defined on $\hat{{\mc S}}_{q^2}$ since ${\bf (W42)}\,\,\gL
L_{q^2}=L_{q^2}\gL$ and $\pp_zL_{q^2}=L_{q^2}\pp_z$.  Moreover the $q^2$ integral
vanishes on the functions from $ker(L_{q^2})$ and is therefore well defined on the
quotient in {\bf (W40)}.  One will now occasionally write in the $q^2$ integral 
an element from $\hat{{\mc S}}_{q^2}$ assuming that it is a representative from the
quotient.  Then if $\hat{\phi}(z)\in\hat{{\mc S}}{q^2}$ one has ${\bf (W43)}\,\,\int
d_{q^2}z\pp_z\hat{\phi}(z)=0$, which follows from {\bf W30)},\,{\bf (W37)}, and
\eqref{11e} via
\bq\label{17e}
\int
d_{q^2}z\pp_z\hat{\phi}(z)=\sum_{-\infty}^{\infty}[\hat{\phi}(q^{2m}-\hat{\phi}(q^{2m+2})-
\hat{\phi}(-q^{2m})+\hat{\phi}(-q^{2m+2})]=
\end{equation}
$$=lim_{M\to\infty}[\hat{\phi}(q^{-2M})-\hat{\phi}(-q^{-2M})]=0$$
From this follows now a $q^2$ integration by parts formula for $k\geq 0$
\bq\label{18e}
\int d_{q^2}z\hat{\phi}\pp_z^k\hat{\psi}(z)=(-1)^kq^{-k(k-1)}\int
d_{q^2}z\pp_z^k\hat{\phi}(z)\hat{\psi}(q^{2k}z)
\end{equation}
\begin{definition}
A $q^2$ distribution f over $\hat{{\mc S}}_{q^2}$ is a continuous linear functional
$f:\,\hat{{\mc S}}_{q^2}\to {\bf C}$ and the space of such distributions is denoted by
$\hat{{\mc S}}'_{q^2}$.  A sequence $f_n\to f\in\hat{{\mc S}}_{q^2}$ if for any
sequence $\phi\in\hat{{\mc S}}_{q^2}$ one has $<f_n,\phi>\to<f,\phi>$.
Distributions determined by
\bq\label{19e}
<f,\phi>=\int_{-\infty}^{\infty} d_{q^2}z\bar{f}(z)\phi(z)=(1-q^2)\sum_{-\infty}^{\infty}
q^{2m}[\bar{f}(q^{2m})\phi(^{2m})+\bar{f}(-q^{2m})\phi(-q^{2m})]
\end{equation}
are called regular.
From \eqref{18e} and {\bf (N18)} one can introduce $q^2$ differentiation in 
$\hat{{\mc S}}'_{q^2}$ via
${\bf (W44)}\,\,<\pp_zf,\phi>=-q^2<\gL f,\pp_z\phi>$.
\end{definition}
\begin{example}
Some examples of $q^2$ distributions are as follows.
\begin{itemize}
\item
\bq\label{20e}
<\gt_{q^2}^{+},\phi>=\int_0^{\infty}d_{q^2}z\hat{\phi}(z)=(1=q^2)
\sum_{-\infty}^{\infty}q^{2m}\phi(q^{2m});
\end{equation}
$$<\gt_{q^2}^{-},\phi>=\int_{-\infty}^0d_{q^2}z\hat{\phi}(z)=(1-q^2)
\sum_{-\infty}^{\infty}q^{2m}\phi(-q^{2m})$$
Thus $\gt_{q^2}^{+}$ and $\gt_{q^2}^{-}$ correspond to functions ${\bf
(W45)}\,\,\gt_{q^2}^{+}(z)=\sum_{-\infty}^{\infty}\hat{\phi}_n^{+}(z)$ and
$\gt_{q^2}^{-}=\sum_{-\infty}^{\infty}\hat{\phi}_n^{-}(z)$.
\item
\bq\label{21e}
<\gd_{q^2},\phi>=\phi(0)=lim_{m\to\infty}\frac{\phi(q^{2m})+\phi(-q^{2m})}{2}
\end{equation}
One can also show easily that ${\bf (W46)}\,\,\pp_z(\gt_{q^2}^{+}(z)-\gt_{q^2}^{-}(z))
=2\gd_{q^2}(z)$.
\item
For arbitrary $k\geq 0$
\bq\label{24e}
<z^{-k-1},\phi>=(-1)^kq^{k(k+1)}\frac{(1-q^2)^k}{(q^2;q^2)_k}<\pp_z^kz^{-1},\phi(z)>=
\end{equation}
$$=\frac{(1-q^2)^k}{(q^2;q^2)_k}\sum_{-\infty}^{\infty}[\pp^k_z\phi(z)|_{z=q^{2m}}-
\pp_z^k\phi(z)|_{z=-q^{2m}}]$$
\item
For arbitrary $\nu>1$
\bq\label{25e}
<z^{\nu}_{+},\phi>=\int_0^{\infty}
d_{q^2}za^{\nu}\hat{\phi}(z)=(1-q^2)\sum_{-\infty}^{\infty}q^{2m(\nu+1)}
\phi(q^{2m})
\end{equation}
and since for any $k\geq 0$ one has ${\bf (W47)}\,\,\pp_z^kz^{\nu}=(-1)^kq^{k(2\nu-k+1)}
[(q^{-2\nu},q^2)_k/(1-q^2)^k]z^{\nu-k}$ one defines
\bq\label{26e}
<z_{+}^{\nu-k},\phi>=\frac{(1-q^2)^{k+1}}{(q^{-2\nu};q^2)_k}\sum_{-\infty}^{\infty}
q^{2m(\nu+1)}\pp_z^k\phi(z)|_{z=q^{2m}}
\end{equation}
\item
Similarly for an arbitrary $\nu>-1$ 
\bq\label{27e}
<z^{\nu}_{-},\phi>=\int_{-\infty}^0
d_{q^2}z(-z)^{\nu}\hat{\phi}(z)=(1-q^2)\sum_{-\infty}^{\infty}q^{2m(\nu+1)}\phi
(-q^{2m})
\end{equation}
and for arbitrary $k\geq 0$
\bq\label{28e}
<z_{-}^{\nu-k},\phi>=(-1)^k\frac{(1-q^2)^{k+1}}
{(q^{-2\mu};q^2)_k}\sum_{-\infty}^{\infty}q^{2m(\nu
+1)}\pp_z^k\phi(z)|_{z=-q^{2m}}
\end{equation}
\end{itemize}
\end{example}
\indent
Let ${\mc S}^{q^2}=\{\psi(s)\}$ generated by a similar class of functions to those in 
${\mc S}_{q^2}$ with $s,\,z$ related as the generators in $\mf{A}\mf{B}$.  Introduce the
same kind of topology via ${\bf (W48)}\,\,|s^k\pp_s^{\ell}\psi(s)|\leq C_{k,\ell}(q)$ for
$k,\ell\geq 0$ and define the map
\bq\label{29e}
{\mc S}_{q^2}\stackrel{L_{q^2}}{\to}\hat{{\mc
S}}_{q^2}\stackrel{\mf{F}_{q^2}}{\to}{\mc S}^{q^2};\,\,\mf{F}_{q^2}\phi(z)=\int
d_{q^2}z\phi(z){}_0\Phi_1(-;0;q^2;i(1-q^2)q^2zs)
\end{equation}
(cf. {\bf (W27)}).  This is the $q^2$ Fourier transform and in the following one sometimes
discards the $L_{q^2}$ notation.  For the inverse and continuity look at the dual space
of skeletons $\hat{{\mc S}}^{q^2}$ and the map
\bq\label{30e}
\mf{F}_{q^2}^{-1}\psi(s)=\frac{1}{2\gT_0}\int E_{q^2}(-i(1-q^2)zs)\psi(s)d_{q^2}s
\end{equation}
where $\mf{F}_{q^2}^{-1}:\,\hat{{\mc S}}^{q^2}\to {\mc S}_{q^2}$.  Consider the diagram
\[\
\begin{CD}
{\mc S}_{q^2}           @>{L_{q^2}}>>      \hat{{\mc S}}_{q^2}\\
@A{\mf{F}^{-1}_{q^2}}AA                     @V{\mf{F}_{q^2}}VV\\
\hat{{\mc S}}^{q^2}     @<{L_{q^2}}<<           {\mc S}^{q^2}      
\end{CD}
\]
where it is shown that 
\bq\label{32e}
\mf{F}_{q^2}L_{q^2}:\,\,{\mc S}_{q^2}\to{\mc S}^{q^2};\,\,\mf{F}_{q^2}^{-1}L_{q^2}:\,\,
{\mc S}^{q^2}\to{\mc S}_{q^2}
\end{equation}
are topological isomorphisms.  Similarly the maps 
\bq\label{33e}
L_{q^2}\mf{F}_{q^2}^{-1}L_{q^2}\mf{F}_{q^2}:\,\,\hat{{\mc S}}_{q^2}\to\hat{{\mc S}}_{q^2}
\end{equation}
and $L_{q^2}\mf{F}_{q^2}L_{q^2}\mf{F}_{q^2}^{-1}$ are identity maps.
To see this one writes
\bq\label{34e}
\int d_{q^2}ze_{q^2}(-i(1-q^2)z){}_0\Phi_1(-;0;q^2;i(1-q^2)q^2zs)=\left\{\begin{array}{cc}
\frac{2}{1-q^2}\gT_0 & s=1\\
0 & s\ne 1
\end{array}\right.
\end{equation}
In accordance with the definition of $q^2$ integrals the integrands must be ordered and
it follows from \eqref{2e}, {\bf (W27)}, and {\bf (W31)} that ${\bf
(W49)}\,\,{}_0\Phi_1(-;0;q^2;i(1-q^2)q^2sz)=\ddg E_{q^2}(i(1-q^2)zs)\ddg$ leading to
${\bf (W50)}\,\,\int d_{q^2}ze_{q^2}(-i(1-q^2)z)\ddg E_{q^2}(i(1-q^2)q^2zs)\ddg$.  Using
then \eqref{2e} one obtains
\bq\label{35e}
\pp_z[e_{q^2}(-i(1-q^2)z)\ddg E_{q^2}(i(1-q^2)zs\ddg]=
\end{equation}
$$=-ie_{q^2}(-i(1-q^2)z)\ddg
E_{q^2}(i(1-q^2)q^2zs)\ddg (1-s)$$
Hence if $s\ne 1$
\bq\label{36e}
\int d_{q^2}ze_{q^2}(-i(1-q^2)z)\ddg E_{q^2}(i(1-q^2)q^2zs)\ddg=
\end{equation}
$$=i(1-s)^{-1}\int d_{q^2}z\pp_x[e_{q^2}(-i(1-q^2)z)\ddg E_{q^2}(i(1-q^2)zs)\ddg]$$
Using \eqref{3e}, \eqref{6e}, and {\bf (W28)} one shows that this vanishes while if 
$s=1$ one gets the value $2(1-q^2)^{-1}\gT_0$ (cf. \eqref{7e} and {\bf (W29)}).  Next one
proves that
\bq\label{37e}
\int E_{q^2}(-i(1-q^2)zs)E_{q^2}(i(1-q^2)q^2s)d_{a^2}s=\left\{\begin{array}{cc}
\frac{2}{1-q^2} & z=1\\
0 & z\ne 1
\end{array}\right.
\end{equation}
Further
\bq\label{38e}
\mf{F}_{q^2}\gL=q^{-2}\gL^{-1}\mf{F}_{q^2};\,\,\mf{F}_{q^2}\pp_z=
-is\mf{F}_{q^2};\,\,\mf{F}_{q^2}z=-iq^{-2}\gL^{-1}\pp_s\mf{F}_{q^2};
\end{equation}
$$\mf{F}_{q^2}^{-1}\gL=q^{-2}\gL^{-1}\mf{F}_{q^2}^{-1};\,\,\mf{F}_{q^2}^{-1}
\pp_s=i\gL^{-1}z\mf{F}_{q^2}^{-1};\,\,\mf{F}_{q^2}^{-1}s=i\pp_z\mf{F}_{q^2}^{-1}$$
This is straightforward using 
\bq\label{39e}
{}_0\Phi_1(-;0;q^2;i(1-q^2)q^2zs)=
\end{equation}
$$=\ddg E_{q^2}(i(1-q^2)q^2zs)\ddg;\,\,E_{q^2}(-i(1-q^2)zs)=\ddg e_{q^2}(i(1-q^2)zs)\ddg$$
\bq\label{40e}
\pp_z\ddg E_{q^2}((1-q^2)azs)\ddg=a\ddg E_{q^2}((1-q^2)aq^2zs)\ddg s;
\end{equation}
$$\pp_s\ddg E_{q^2}((1-q^2)azs)\ddg=az\ddg E_{q^2}((1-q^2)aq^2zs)\ddg;$$
$$\pp_z\ddg e_{q^2}((1-q^2)azs)\ddg=a\ddg e_{q^2}((1-q^2)azs)\ddg s;$$
$$\pp_s\ddg e_{q^2}((1-q^2)azs)\ddg =az\ddg e_{q^2}((1-q^2)azs)\ddg$$
\\[3mm]\indent
Finally one can prove \eqref{32e}-\eqref{33e}.  First from \eqref{38e} one has
\bq\label{41e}
\mf{F}_{q^2}z^k\pp_z^{\ell}\phi(z)=(-i)^{k+\ell}q^{-2k}(\gL^{-1}\pp_s)^ks^{\ell}
\mf{F}_{q^2}\phi(z)
\end{equation} 
and on the other hand
\bq\label{42e}
\pp_s^ks^{\ell}=(-1)^kq^{k(2\ell
-k+1)}\sum_{j=0}^k(-1)^jq^{j(j-1)}\frac{(q^{-2\ell};q^2)_j}
{(1-q^2)^{k-j}}\left[\begin{array}{c}
k\\
j\end{array}\right]_{q^2}s^{\ell-k+j}\pp_s^j
\end{equation}
Hence if $\phi(z)$ satisfies {\bf (W37)} then its image $\mf{F}_{q^2}\phi(z)$ satisfies
{\bf (W48)}; this means that the image of convergent sequences in $\hat{{\mc S}}_{q^2}$
converges in $\hat{{\mc S}}^{q^2}$ and a similar statement for $\mf{F}_{q^2}^{-1}$ is
proved in the same manner using \eqref{38e}.  Next consider the action of the Fourier
operators on the basis
\bq\label{43e}
\mf{F}_{q^2}\ci\mf{F}_{q^2}^{-1}\hat{\psi}_n^{\pm}(x)=\hat{\psi}_n^{\pm}(s);\,\,
\mf{F}_{q^2}^{-1}\ci\mf{F}_{q^2}\hat{\phi}_n^{\pm}(z)=\hat{\phi}_n^{\pm}(z)
\end{equation}
Take e.g. 
$$\mf{F}_{q^2}\ci\mf{F}_{q^2}^{-1}\hat{\psi}_n^{+}(s)=\frac{1}{2\gT_0}\int
d_{q^2}z\left(\int E_{q^2}(-i(1-q^2)z\xi)\hat{\psi}_n^{+}(\xi)d_{q^2}\xi\right){}_0\Phi_1
(-;0;q^2;i(1-q^2)zs)=$$
\bq\label{44e}
=\frac{1-q^2}{2\gT_0}q^{2n}\int d_{q^2}ze_{q^2}(-i(1-q^2)q^{2n}z){}_0\Phi_1(-;0;q^2;i
(1-q^2)q^2zs)=
\end{equation}
$$=\frac{1-q^2}{2\gT_0}\int d_{q^2}ze_{q^2}(-i(1-q^2)z){}_0\Phi_1(-;0;q^2;i(1-q^2)
q^{-2n+2}zs)$$
It follows from \eqref{34e} that 
\bq\label{45e}
\mf{F}_{q^2}\ci\mf{F}_{q^2}^{-1}\hat{\psi}_n^{+}(s)=\left\{\begin{array}{cc}
1 & s=q^{2n}\\
0 & x\ne q^{2n}
\end{array}\right.
\end{equation}
and one arrives at \eqref{43e} (using \eqref{37e} for the second equation).
\begin{definition}
The $q^2$ Fourier transform of a $q^2$ distribution $f\in\hat{{\mc S}}_{q^2}'$ is the 
$q^2$ distribution $g\in (\hat{{\mc S}}^{q^2})'$ defined via $<g,\psi>=<f,\phi>$ where
$\phi\in {\mc S}_{q^2}$ is arbitrary and $\psi\in\hat{{\mc S}}^{q^2}$ is its $q^2$
Fourier transform.
\end{definition}
\indent
Suppose that $zf(z)$ is absolutely $q^2$ integrable for the $q^2$ distribution $f$ and
let $\phi(z)=\mf{F}_{q^2}^{-1}\hat{\psi}(s)$; then
\bq\label{46e}
<f,\phi>=\frac{1}{2\gT_0}\int d_{q^2}z\bar{f}(z)\int E_{q^2}(-i(1-q^2)zs)\psi(s)d_{q^2}s=
\end{equation}
$$=\frac{1}{2\gT_0}\int\ol{\int d_{q^2}zf(z)E_{q^2}(i(1-q^2)zs)}\psi(s)d_{q^2}s=<g,\psi>$$
This means that the $q^2$ distribution $g\sim 
g(s)=(1/2\gT_0)\int d_{q^2}zf(z)E_{q^2}(i(1-q^2)zs)$.  Similar calculations lead to the
following relations in the space of $q^2$ distributions
\bq\label{48e}
\mf{F}'_{q^2}\gL=q^{-2}\gL^{-1}\mf{F}'_{q^2};\,\,\mf{F}'_{q^2}\pp_z=
-i\gL^{-1}s\mf{F}'_{q^2};\,\,\mf{F}'_{q^2}z=-i\pp_s\mf{F}'_{q^2};
\end{equation}
$$(\mf{F}_{q^2}')^{-1}\gL=q^{-2}\gL^{-1}(\mf{F}'_{q^2})^{-1};\,\,
(\mf{F}'_{q^2})^{-1}\pp_s=iz(\mf{F}'_{q^2})^{-1};\,\,(\mf{F}'_{q^2})^{-1}s
=iq^{-2}\gL^{-1}\pp_z(\mf{F}'_{q^2})^{-1}$$
\begin{example}
As examples consider
\begin{itemize}
\item
From {\bf (W35)} follows
\bq\label{49e}
\mf{F}'_{q^2}z^{-1}=i\,sgn(s)=i(\gt_{q^2}^{+}-\gt_{q^2}^{-})
\end{equation}
\item
From \eqref{48e}, \eqref{49e}, and {\bf (W46)} one obtains
\bq\label{50e}
\mf{F}'_{q^2}1=\mf{F}'_{q^2}zz^{-1}=-i\pp_s\mf{F}'_{q^2}z^{-1}=\pp_s(\gt_{q^2}^{+}-
\gt_{q^2}^{-})=2\gd_{q^2}
\end{equation}
\item
\bq\label{51e}
\mf{F}'_{q^2}(\gt_{q^2}^{+}-\gt_{q^2}^{-})=\frac{i(1-q^2)}{\gT_0}
\sum_{-\infty}^{\infty}q^{2m}sin_{q^2}(91-q^2)q^{2m}s)=\frac{is^{-1}}{\gT_0}
\end{equation}
\item
From {\bf (W46)}, \eqref{48e}, and \eqref{51e} results
\bq\label{52e}
\mf{F}'_{q^2}\gd=\frac{1}{2}\mf{F}'_{q^2}\pp_z(\gt_{q^2}^{+}-\gt_{q^2}^{-})=
-\frac{i}{2}\gL^{-1}s\mf{F}'_{q^2}(\gt_{q^2}^{+}=\gt_{q^2}^{-})=\frac{1}{2\gT_0}
\end{equation}
\item
From \eqref{20e} and \eqref{51e} one obtains 
\bq\label{53e}
\mf{F}'_{q^2}\gt_{q^2}^{+}=\frac{1}{2}\mf{F}'_{q^2}(\gt_{q^2}^{+}-
\gt_{q^2}^{-}+1)=\frac{is^{-1}}{2\gT_0}+\gd_{q^2};
\end{equation}
$$\mf{F}'_{q^2}\gt_{q^2}^{-}=\frac{1}{2}\mf{F}'_{q^2}(-\gt^{+}_{q^2}+
\gt_{q^2}^{-}+1)=-\frac{is^{-1}}{2\gT_0}+\gd_{q^2}$$
\end{itemize}
\end{example}
\indent
Further for $n\geq 0$ arbitrary
\bq\label{54e}
\mf{F}'_{q^2}z^n=2i^nq^{-n(n+1)}\frac{(q^2;q^2)_n}{(1-q^2)^n}s^{-n}\gd_{q^2}(s);
\end{equation}
$$\mf{F}'_{q^2}z^{-n-1}=i^{n+1}\frac{(q^2;q^2)_n}{(1-q^2)^n}s^nsgn(s)$$
and (cf. \cite{oe} for proofs)
\bq\label{55e}
\mf{F}'_{q^2}z_{+}^{\nu-1}=\frac{e_{q^2}(q^2)E_{q^2}(-q^{2(1-\nu)})}
{2\gT_0}(\bar{c}_{\nu}s_{-}^{-\nu}+c_{\nu}s_{+}^{-\nu});
\end{equation}
$$\mf{F}'_{q^2}z_{-}^{\nu-1}=-\frac{e_{q^2}(q^2)E_{q^2}(-q^{2(1-\nu)})}
{2\gT_0}(c_{\nu}s_{-}^{-\nu}+\bar{c}_{\nu}s_{+}^{-\nu})$$
$$c_{\nu}=\sum_{-\infty}^{\infty}\frac{q^{2\nu
m}(q^{-2m}+i(1-q^2))}{(1-q^2)^{-1}q^{-2m}+(1-q^2)q^{2m}}$$
Further examples and calculations can be found in \cite{oe}.

\section{SOME PRELIMINARY CALCULATIONS}
\renewcommand{\theequation}{6.\arabic{equation}}
\setcounter{equation}{0}

We go back now to Section 3 and try to construct counterparts for some of the
classical ideas there.  
\subsection{{\bf Background machinery.}}
First however we
consider $SL_q(2)$ with relations (3)-(7) and \eqref{71a}.
$$x^1x^2=qx^2x^1;\,\,\pp_1x^1=1+q^2x^1\pp_1+(q^2-1)x^2\pp_2;
\,\,\pp_1x^2=qx^2\pp_1;\,\,\pp_2x^1=qx^1\pp_2;$$
$$\pp_2x^2=1+q^2x^2\pp_2;\,\,
\pp_1\pp_2=q^{-1}\pp_2\pp_1;\,\,dx^1x^1=q^2x^1dx^1;$$
\bq\label{1b}
dx^1x^2=qx^2dx^1+(q^2-1)x^1dx^2;\,\,
dx^2x^1=qx^1dx^2;\,\,dx^2x^2=q^2x^2dx^2
\end{equation}
Here $x^1$ and $x^2$ generate a quantum plane V or $H=SL_q(2)$ comodule where
$SL_q(2)$ is described via the R matrix \eqref{45a} and relations \eqref{44a} (with
$det_qT=1$); we will use H and $\mf{A}=\mf{O}(SL_q(2))$ interchangeably at times
(see below for more on this).
The coaction is given by {\bf (Z53)} in the form $\go(x^i)=\sum
T^i_k\ot x^k\sim\gD_V(x^i)\in H\ot V$.  Note the $\pp_i$ transform covariantly via
$\go(\pp_i)=\sum\hat{S}^{\ell}_i\ot\pp_{\ell}$ with 
$\hat{S}=\widetilde{\tl{T}^{-1}}$ where $\tl{T}^a_b=T^b_a$; here $\tl{T}^{-1}\ne
\widetilde{T^{-1}}$ since for $SL_q(2)$ one has
$$T^{-1}=\left(\begin{array}{cc}
d & -q^{-1}b\\
-qc & a
\end{array}\right);\,\,\tl{T}=\left(\begin{array}{cc}
a & c\\
b & d
\end{array}\right);\,\,
\tl{T}^{-1}=\left(\begin{array}{cc}
d & -q^{-1}c\\
-qb & a
\end{array}\right);$$
\bq\label{2b}
\widetilde{T^{-1}}=\left(\begin{array}{cc}
d & -qc\\
-q^{-1}b & a
\end{array}\right);\,\,\widetilde{\tl{T}^{-1}}=\left(\begin{array}{cc}
d & -qb\\
-q^{-1}c & a
\end{array}\right)
\end{equation}
We could think of the $\pp_i$ also generating an H comodule W with
$\gD_W(\pp_i)=\go(\pp_i)\in H\ot W$.  Recall also from {\bf (Z49)} that $\gD
T^j_{\ell}=\sum T^j_r\ot T^r_{\ell}$ and following \cite{ka}, p. 97, one expects
$$\gD(a)=a\ot a+b\ot c;\,\,\gD(b)=a\ot b+b\ot d;$$
\bq\label{3b}
\gD(c)=c\ot a+d\ot c;\,\,
\gD(d)=c\ot b+d\ot d
\end{equation}
(e.g. $T^1_2=b$ so $\gD(b)=T^1_1\ot T^1_2+T^1_2\ot T^2_1=a\ot b+b\ot d$ as
indicated).  We note also the agreement of \eqref{1b} or more generally (3)-(5) with
the constructions in \cite{ka}, pp. 468-469 for $\mf{A}=\mf{O}(SL_q(n))$ (for $n=2$
this is the Hopf algebra $\mf{A}$ generated by $1,a,b,c,d$ modulo the relations
\eqref{44a} and $det_qT=1$) and we refer to it as the quantum algebra $SL_q(n)$; more
correctly it is the coordinate algebra of $SL_q(n)$.  We recall $\gep(T)=id$ and
$S(T)=T^{-1}$ which translates into
\bq\label{4b}
S(a)=d;\,\,S(b)=-q^{-1}b;\,\,S(c)=-qc;\,\,S(d)=a
\end{equation}
As a left comodule the quantum plane $V\sim \mf{O}({\bf C}_q^2)$ obeys
$\go(x^i)=\gD_V(x^i)=\phi_L(x^i)=\sum T^i_k\ot x^k$ which can be written out as
\bq\label{5b}
\phi_L(x^1)=a\ot x^1+ b\ot x^2;\,\,\phi_L(x^2)=c\ot x^1+d\ot x^2
\end{equation}
V is also a right comodule of $\mf{A}=H$ via
\bq\label{5cc}
\phi_R(x^1)=x^1\ot a+x^2\ot c;\,\,\phi_R(x^2)=x^1\ot b+x^2\ot d
\end{equation}
For completeness we specify the $\gD_W(\pp_i)$ above via
\bq\label{6b}
\gD_W(\pp_1)=\hat{S}_1^1\ot\pp_i+\hat{S}^2_1\ot\pp_2=d\ot\pp_1-qb\ot\pp_2;
\end{equation}
$$\gD_W(\pp_2)=\hat{S}^1_2\ot\pp_1+\hat{S}_2^2\ot\pp_2=-q^{-1}c\ot\pp_1+a\ot\pp_2$$
\indent
We go next to \cite{ka}, pp. 468-469, where particular covariant FODC $\gG_{\pm}$ on 
$\mf{O}({\bf C}^n_q)$ are discussed.  The algebra $\mf{O}({\bf C}^n_q)$ is defined as
in \eqref{1b} with $\gG_{+}$ described via 
\bq\label{7b}
x_i\cdot dx_j=qdx_j\cdot x_i+(q^2-1)dx_i\cdot x_j\,\,(i<j);\,\,x_i\cdot dx_i=q^2
dx_i\cdot x_i;
\end{equation}
$$x_j\cdot dx_i=qdx_i\cdot x_j\,\,(i<j);\,\,dx_i\wg dx_j=-q^{-1}dx_j\wg
dx_i\,\,(i<j);\,\,dx_i\wg dx_i=0$$
\bq\label{8b}
x_ix_j=qx_jx_i\,\,(i<j);\,\,\pp_i\pp_j=q^{-1}\pp_j\pp_i\,\,(i<j);\,\,\pp_ix_j=qx_j\pp_i
\,\,(i\ne j);
\end{equation}
$$\pp_ix_i-q^2x_i\pp_i=1+(q^2-1)\sum_{j>i} x_j\pp_j$$
and $\gG_{-}$ arises by replacing $q$ by $q^{-1}$ and $i<j$ by $j<i$ in the formulas
\eqref{7b}-\eqref{8b}.  Note for $q\ne 1$, $\gG_{+}$ and $\gG_{-}$ are not isomorphic
and for $q=1$ they both give the ordinary differential calculus on the correponding
polynomial algebra ${\bf C}[x_1,\cdots,x_n]$.
Let us summarize some facts about $\gG_{\pm}$ in
\begin{example}
There are two distinguished FODC, $\gG_{\pm}$ on $\mf{O}({\bf C}_q^2)$. 
For both calculi the set of differentials $\{dx,dy\}$ is a basis for the right (and
the left) $\mf{O}({\bf C}_q^2)$ module of first order forms.  Hence for any
$z\in\mf{O}({\bf C}_q^2)$ there exist uniquely determined elements $\pp_x(z)$ and
$\pp_y(z)$, called partial derivatives, such that ${\bf
(W51)}\,\,dz=dx\cdot\pp_x(z)+dy\cdot\pp_y(z)$.  The bimodule structures of $\gG_{\pm}$
are given via
\bq\label{35b}
\gG_{+}:\,\, xdy=qdy\cdot x+(q^2-1)dx\cdot y;\,\,ydx=qdx\cdot y;\,\,xdx=q^2dx\cdot
x;\,\,ydy=q^2dy\cdot y
\end{equation}
$$\gG_{-};\,\,ydx=q^{-1}dx\cdot y+(q^{-2}-1)dy\cdot x;\,\,xdy=q^{-1}dy\cdot
x;\,\,xdx=q^{-2}dx\cdot x;\,\,ydy=q^{-2}dy\cdot y$$
From this one sees that $\eta_{+}=y^{-2}xdx$ and $\eta_{-}=x^{-2}ydy$ are nonzero
central elements of the bimodules $\gG_{+}$ and $\gG_{-}$ respectively (recall central
means $\eta z=z\eta$ for all z).  One notes also that the relations for $\gG_{+}$ go
into those for $\gG_{-}$ if we interchange the coordinates x and y and $q\to q^{-1}$. 
The partial derivatives $\pp_x$ and $\pp_y$, considered as linear mappings of 
$\mf{O}({\bf C}_q^2)$, and the coordinate functions $x,\,y$, acting on $\mf{O}({\bf
C}_q^2)$ by left multiplication, satisfy the relations
\bq\label{36b}
\gG_{+}:\,\,\pp_xy=qy\pp_x;\,\,\pp_yx=qx\pp_y;\,\,\pp_xx-q^2x\pp_x=
1+(q^2-1)y\pp_y;\,\,\pp_yy-q^2y\pp_y=1;
\end{equation}
$$\gG_{-}:\,\,\pp_xy=q^{-2}y\pp_x;\,\,\pp_yx=q^{-1}x\pp_y;\,\,\pp_x-q^{-2}x\pp_x=1;
\,\,\pp_yy-q^{-2}y\pp_y=1+(q^{-2}-1)x\pp_x$$
From these formulas one derives by induction the expressions for the actions of
$\pp_x$ and $\pp_y$ on general elements of $\mf{O}({\bf C}_q^2)$ and for polynomials
$g$ and $h$ one has
\bq\label{37b}
\gG_{+}:\,\,\pp_x(g(y)h(x))=g(qy)D_{q^{-2}}(h)(x);\,\,\pp_y(g(y)h(x))=
D_{q^{-2}}(g)(y)h(x);
\end{equation}
$$\gG_{-}:\,\,\pp_x(g(x)h(y))=D_{q^{-2}}(g)(x)h(y);\,\,\pp_y(g(x)h(y))=
g(q^{-1}x)D_{q^{-2}}(h)(y)$$
We note here that in \cite{sz} one points out that rules for $\pp_i$ (or 
\eqref{1b}) actually can be regarded as q-deformed Leibnitz rules (and also
in part as q-deformed Heisenberg relations).  Thus consider 
in $\gG_{-}$ for example $\pp_xx-q^{-2}x\pp_x=1$ and replace $q^{-2}$ by $q$ for this
illustration.  Then $\pp_xx-qx\pp_x=1$ and with $p=-i\pp_x$ this becomes ${\bf
(W52)}\,\,px-qxp=-i$.  Then for $p=-i\bar{p}$ one requires $\bar{x}p-qp\bar{x}=i$ or 
$\bar{x}p-(1/q)p\bar{x}=-(1/q)$ which involves the introduction of a new element
$\bar{x}$ into the algebra (cf. here \eqref{79a} where $\gL$ was introduced).  Note
as a Leibnitz rule $\pp_xx-qx\pp_x=1$ can be related to \eqref{9b}-\eqref{10b} for
example. Let us work from $\gG_{+}$ now and then $\pp_x=q^2x\pp_x+1+(q^2-1)y\pp_y$. 
First to derive $\pp_xx^n=[[n]]_{q^2}x^{n-1}$ we have
\bq\label{38b}
\pp_xx=1;\,\,\pp_xx^2=q^2x\pp_xx+x=(q^2+1)x;\,\,\pp_xx^3=q^2x\pp_xx^2+x^2=
\end{equation}
$$=(q^4+q^2+1)x^2=\frac{q^6-1}{q^2-1}x^2;\,\,\cdots$$
Consequently for a polynomial $f(x)=\sum_0^N a_nx^n$ one has
\bq\label{39b}
\pp_xf(x)=\sum_0^Na_n[[n]]_{q^2}x^{n-1}=f'_{q^{2}}(x)\sim D_{q^2}f(x)
\end{equation}
as in \eqref{19b}.  For Leibnitz a simple calculation gives for polynomials $f,\,g$
\bq\label{40b}
\pp_x(fg)=D_{q^2}(fg)=\frac{(fg)(x)-(fg)(q^2x)}{(1-q^2)x}=
\end{equation}
$$=\frac{g(x)[f(x)-f(q^2x)]+f(q^2x)[g(x)-g(q^2x)]}{(1-q^2)x}=g(x)\pp_xf(x)+
f(q^2x)\pp_xg(x)$$
\end{example}
\begin{definition}
The q-Weyl algebra ${\mc A}_q(n)$ is the unital algebra with $2n$ generators
$x_1,\cdots,x_n,\pp_1,\cdots,\pp_n$ determined by the relations \eqref{8b}.  $\mf{D}$
is the unital algebra generated by $\pp_1,\cdots,\pp_n$ with
$\pp_i\pp_j=q^{-1}\pp_j\pp_i$ for $i<j$.  One can consider $\mf{A}=\mf{O}({\bf
C}^n_q)$ and $\mf{D}$ as subalgebras of ${\mc A}_q(n)$ and the set of monomials
$\{x_1^{k_1}\cdots x_n^{k_n}\pp_i^{m_1}\cdots\pp_n^{m_n}\}$ as a basis.
The element $D=\sum x_i\pp_i\in {\mc
A}_q(n)$ is called the Euler derivation.  
\end{definition}
Before going further we want to clarify the role of $U_q(s\ell_2)$ in terms of
differential operators (cf. \cite{fa,kc,wzz}).
\begin{definition}
Following \cite{ka} $U_q(s\ell_2)$ is defined as the unital associative algebra over 
{\bf C} with generators $E,\,F,\,K,\,K^{-1},$ satisfying 
\bq\label{20b}
KK^{-1}=K^{-1}K=1;\,\,KEK^{-1}=q^2E;\,\,KFK^{-1}=q^{-2}F;\,\,[E,F]=\frac{K-K^{-1}}{q-q^{-1}}
\end{equation}
One can take $\{F^{\ell}K^mE^n;\,\,m\in{\bf Z},\,\,\ell,n\in{\bf N}_0\}$ or 
$\{E^nK^mF^{\ell};\,\,m\in{\bf Z},\,\,\ell,n\in{\bf N}_0\}$ as a vector space basis of
$U_q(s\ell_2)$.  The quantum Casimir element 
\bq\label{21b}
C_q=EF+\frac{Kq^{-1}+K^{-1}q}{(q-q^{-1})^2}=FE+\frac{Kq+K^{-1}q^{-1}}{(q-q^{-1})^2}
\end{equation}
lies in the center of $U_q(s\ell_2)$ and if q is not a root of unity $C_q$ generates
the center.
\end{definition}
\indent
There is a unique Hopf algebra structure on $U_q(s\ell_2)$ built upon
$$\gD(E)=E\ot K+1\ot E;\,\,\gD(F)=F\ot 1+K^{-1}\ot F;\,\,\gD(K)=K\ot
K;\,\,S(K)=K^{-1};$$
\bq\label{22b}
S(E)=-EK^{-1};\,\,S(F)=-KF;\,\,\gep(K)=1;\,\,\gep(E)=\gep(F)=0
\end{equation}
\indent
One recalls that $U(s\ell_2)=U(\mf{g})$ is the tensor algebra $T(\mf{g})$ modulo the
two sided ideal generated by $x\ot y-y\ot x-[x,y]$ for $x,y\in\mf{g}$ and here
$U(s\ell_2)$ is generated by elements $x_{\pm}$ and $y$ with
\bq\label{25b}
[y,x_{+}]=2x_{+};\,\,[y,x_{-}]=-2x_{-};\,\,[x_{+},x_{-}]=y;
\end{equation}
$$\gD(x)=x\ot 1+1\ot
x;\,\,\gep(x)=0;\,\,S(x)=-x$$
Note that for $q=exp(h)$ and $K=exp(hH)$, in the limit $h\to 0$ \eqref{20b} implies
${\bf (W53)}\,\,[H,E]=2E,\,\,[H,F]=-2F,\,\,[E,F]=H$ corresponding to $E\sim x_{+},\,\,
F\sim x_{-},\,\,y\sim H$.  However this is not quite adequate for comparing
$U_q(s\ell_2)$ to $U(s\ell_2)$.  One considers instead $\tl{U}_q(s\ell_2)$ with
generators $E,\,F,\,K,\,K^{-1}$ and ${\bf (W54)}\,\,G=(q-q^{-1})^{-1}(K-K^{-1})$
satisfying
\bq\label{26b}
[G,E]=E(qK+q^{-1}K^{-1});\,\,[G,F]=-(qK+q^{-1}K^{-1});\,\,[E,F]=G
\end{equation}
Then or $q^2\ne 1,\,\,\tl{U}_q(s\ell_2)$ is isomorphic to $U_q(s\ell_2)$ via $E\to
E,\,F\to F,\,K\to K,$ and $G\to (q-q^{-1})(K-K^{-1})$ (note $\gD(G)=G\ot K+K^{-1}\ot
G,\,\,\gep(G)=0$, and $S(G)=-G$).  Further these formulas hold at $q=\pm 1$ for
$\tl{U}_q$ (which is excluded for $U_q$) and $\tl{U}_1(s\ell_2)$ is considered as the
classical limit of $\tl{U}_q(s\ell_2)$.  In fact ${\bf
(W55)}\,\,\tl{U}_1(s\ell_2)\simeq U(s\ell_2)\ot{\bf C}{\bf Z}_2$ (see \cite{ka}
for more details).
\\[3mm]\indent
It will be useful to recall here some relations involving q-special functions and their
origins based on q-groups and q-algebras (see especially \cite{fw,kd,kf,ki}).
We will illustrate matters using $s\ell_q(2)$ and $SL_q(2)$ and only deal with simple
situations.  First recall
\bq\label{41b}
(a;q)_{\ga}=\frac{(a;q)_{\infty}}{(aq^{\ga};q)_{\infty}};\,\,(a;q)_{\infty}=
\prod_0^{\infty}(1-aq^k)\,\,|q|<1
\end{equation}
where $(a;q)_n=(1-a)(1-aq)\cdots (1-aq^{n-1})$.  There are identities
\bq\label{42b}
q^{(1/2)n(n-1)}(a^{-1}q^{1-n};q)_n=(-a^{-1})^n(a;q)_n;\,\,\sum_0^{\infty}\frac{(a;q)_n}
{(q;q)_n}z^n=\frac{(az;q)_{\infty}}{(z;q)_{\infty}}
\end{equation}
where $|z|<1$ and $|q|<1$.  Note also that the binomial symbol in \eqref{13b}
can be written as
\bq\label{422b}
\nm_q=\frac{(q;q)_n}{(q;q)_m(q;q)_{n-m}}
\end{equation}
There are two q-exponential functions, namely
\bq\label{43b}
e_q(z)=\sum_0^{\infty}\frac{z^n}{(q;q)_n}=\frac{1}{(z;q)_{\infty}};\,\,E_q(z)=
\sum_0^{\infty}\frac{q^{(1/2)n(n-1)}z^n}{(q;q)_n}=(-z;q)_{\infty}
\end{equation}
One notes that $(\bl)\,\,e_q(z)E_q(-z)=1$ and as $q\to 1^{-},\,\,e_q(z(1-q))\to
exp(z)$ and $E_q(z(1-q))\to exp(z)$.  Define ${\bf (W56)}\,\,T_qf(z)=f(qz)$ and
$D^{+}_z=z^{-1}(1-T_q)$ with $D_z^{-}=z^{-1}(1-T^{-1}_q)$.  Then $(1-q)^{-1}D^{+}_z$ 
and $(1-q^{-1})^{-1}D^{-}_z\to d/dz\,\,(q\to 1)$ while ${\bf (W57)}\,\,
D^{+}e_q(z)=e_q(z)$ and $D^{-}_zE_q(z)=-q^{-1}E_q(z)$.  The basic hypergeometric
function is defined via \eqref{44b} in Section 2.
\\[3mm]\indent
Now there are two points of view.  One can work directly with representations of 
$U_q(s\ell_2))\sim U_q(\mf{g})$ or with representations of the coordinate algebra
$\mf{A}(SL_q(2))\sim\mf{A}(G_q)$ which can be considered as a subalgebra of 
$U_q(\mf{g})'$.  There are representations in terms of first or second order 
q-difference operators involving functions of a complex variable z and we mention first
(cf. \cite{fw}) that $s\ell_q(2)$ is generated by 
\bq\label{45b}
J_{+}J_{-}-q^{-1}J_{-}J_{+}=\frac{1-q^{2J_3}}{1-q};\,\,[J_3,J_{\pm}]=\pm J_{\pm}
\end{equation} 
and redefining ${\bf (W58)}\,\,\tl{J}_{\pm}=q^{-(1/2)(J_3\mp (1/2))},\,\,\tl{J}_3=J_3$
this can be written as 
\bq\label{46b}
[\tl{J}_{+},\tl{J}_{-}]=\frac{q^{2\tl{J}_3}-q^{-2\tl{J}_3}}{q-q^{-1}};\,\,
[\tl{J}_3,\tl{J}_{\pm}]=\pm\tl{J}_{\pm}
\end{equation}
Then e.g. the representation $D_q(u,m_0)$ is characterized by two complex constants
$u$ and $m_0$ such that neither $m_0+u$ nor $m_0-u$ is an integer, and $\Re m_0<1$.
On the space $\mf{H}$ of finite linear combinations of $z^n\,\,(n\in{\bf Z})$ the
generators are realized as
\bq\label{47b}
J_{+}=q^{(1/2)(m_0-u+1)}\left[\frac{z^2}{1-q}D^{+}_z-\frac{(1-q^{u-m_0})z}{1-q}\right];
\end{equation}
$$J_{-}=-q^{(1/2)(m_o-u+1)}\left[\frac{1}{1-q}D^{+}_z+\frac{1-q^{u+m_0}}{(1-q)z}T_q
\right];\,\,J_3=m_0+z\frac{d}{dz}$$
The basis vectors $f_m$ in $\mf{H}$ are still defined by $f_m=z^n$ for $m=m_0+n$ and
all integers n and
\bq\label{48b}
J_{+}f_m=q^{(1/2)(u-m_0+1}\frac{1-q^{m-u}}{1-q}f_{m+1};
\end{equation}
$$J_{-}f_m=-q^{(1/2)(m_0-u+1)}\frac{1-q^{m+u}}{1-q}f_{m-1};\,\,J_3f_m=mf_m$$

\subsection{{\bf Q-difference equations.}}

A number of q-difference equations have been treated in the context of symmetries (see 
e.g. \cite{bi,bm,cah,fw,lz,lzz,nb}) and we mention a few results.  Note first that we are
primarily concerned with wave equations in terms of determining transmutation operators
and the wave equation has a number of peculiarities (cf. \cite{fw})  Thus recall {\bf
(W56)} and note
\bq\label{q2}
D_z^{+}e_q(\gl z)=\gl e_q(\gl z);\,\,D_z^{+}E_q(\gl z)=\gl E_q(\gl z);
\end{equation}
$$D_z^{-}e_q(\gl z)=-q^{-1}\gl e_q(q^{-1}\gl z);\,\,D_z^{-}E_q(\gl z)=-q^{-1}\gl E_q(\gl
z)$$
Set also ${\bf (W59)}\,\,{\mc D}_q^{+}=z^{-1}(1-T_q^2)$ and ${\mc
D}_q^{-}=z^{-1}(1-T_q^{-2})$ (so $[1/(1-q^{\pm 2})]{\mc D}_z^{\pm}\to d/dz$ as $q\to
1^{-}$.  Now consider a wave equation in light cone coordinates ${\bf (W60)}\,\,
\pp_1\pp_2\phi=0$.  This has an infinite dimensional symmetry algebra generated by $v^0_m=
x_1^m\pp_1$ and $w^0_m=x_2^m\pp_2$ for ($m\in{\bf Z}$).  In fact the whole
$W_{1+\infty}\oplus W_{1+\infty}$ algebra generated by $V_m^k=x_1^m\pp_1^{k+1}$ and
$w_m^k=x_2^m\pp_2^{k+1}$ for $k\in{\bf Z}_{+}$ maps solutions of {\bf (W60)} into
solutions (definition of a symmetry) and $W_{1+\infty}$ without center corresponds to
$U({\mc E}(2))$.  For the q-difference version ${\bf
(W61)}\,\,D_1^{+}D_2^{+}\phi(x_1,x_2)=0$ the elements $V_m^k=x_1^m(D_1^{+})^{k+1}$ and
$W_m^k=x_2^m(D_2^{+})^{k+1}$ map solutions into solutions and each set $V_m^k$ or $W_m^k$
generates a q-deformation of $W_{1+\infty}$. However for the equation ${\bf
(W62)}\,\,[(D_t^{+})^2-(D_x^{+})^2]\phi(t,x)=0$ the situation is quite different.  There
is still an infinite set of symmetry operators involving polynomials or arbitrary degree
in $t$ and $x$ times powers of $D_t^{+}$ and
$D^{+}_x$ but a general expression seems elusive.  This is in contrast to solving wave
equations $(\pp_t^2-\pp_x^2)\phi=0$ where there is conformal invariance in $t+x$ and
$t-x$.  One notes that $t,x)\to (t+x,t-x)$ does not preserve the exponential
2-dimensional lattice and light cone coordinates seem more appropriate for q-difference
wave equations.
\\[3mm]\indent
In any event (going back to Definition 5.1)
if $f_n\in \mf{A}=\mf{O}({\bf C}^n_q)$ and $g_n\in \mf{D}$ are homogeneous of
degree n then
\bq\label{9b}
Df_n=[[n]]_{q^2}f_n+q^{2n}f_nD;\,\,g_nD=[[n]]_{q^2}g_n+q^{2n}Dg_n
\end{equation}
where $[[n]]_{q^2}=(q^2-1)^{-1}(q^{2n}-1)$.  The partial derivatives $\pp_i$ of 
$\gG_{+}$ act on $\mf{A}$ by the rule
\bq\label{10b}
\pp_i(f_n(x_n)\cdots f_1(x_1))=f_n(qx_n)\cdots f_{i+1}(qx_{i+1})D_{q^2}f_i(x_i)
f_{i-1}(x_{i-1})\cdots f_1(x_1)
\end{equation}
where the $f_i$ are polynomials in one variable and $D_{q^2}$ is the $q^2$ derivative
\bq\label{11b}
D_{q^2}f(x)=\frac{f(x)-f(q^2x)}{x-q^2x}=\frac{f(q^2x)-f(x)}{(q^2-1)x}
\end{equation}
This would seem to be enough now for making calculations.  We will look first at
possible quantum analogues of the wave equation $(\pp_x^2-\pp_y^2)\phi=0$ and quantum
versions of solutions $\phi=F(x-y)+G(x+y)$ (alternatively $(\pp_{xy}+\pp_{yx})\phi=0$ with
solutions
$\phi=\psi(x)+\chi(y)$).  Now the rule
\eqref{10b} is apparently saying that when calculating the action of $\pp_i\in\gG_{+}$ on
functions in $\mf{A}$ one will drop expressions with $\pp_j$ on the right since e.g.
$\pp_j(1)=0$.  Then e.g. using \eqref{1b}, $\pp_1(x_2)=qx_2\pp_1=0$ and $\pp_1(x_1)
=1+q^2x_1\pp_1=1$ while according to e.g. \cite{ka} one has $\pp_1(x_1^k)=[[k]]_{q^2}
x_1^{k-1}=(q^2-1)^{-1}(q^{2k}-1)x_1^{k-1}$.  Note here also
\bq\label{12b}
D_{q^2}x_i^k=\frac{x_1^k-q^{2k}x_1^k}{(1-q^2)x_1}=\frac{q^{2k}-1}{q^2-1}x_1^{k-1}=
[[k]]_{q^2}x_1^{k-1}
\end{equation}
and this is consistent with $x_1\pp_1x_1^k=[[k]]_{q^2}x_1^k$ from \eqref{9b}.  In
order to compute derivatives of functions like $(\ga x_1+\gb x_2)^n$ we recall the 
q-binomial theorem
\bq\label{13b}
(v+w)^n=\sum_0^n\left[\begin{array}{c}
n\\
m\end{array}\right]_{q^2}w^mv^{n-m}=\sum_0^n\left[\begin{array}{c}
n\\
m\end{array}\right]_{q^{-2}}v^mw^{n-m}
\end{equation}
$$ \left[\begin{array}{c}
n\\
m\end{array}\right]_{q^2}=\frac{[n]_q!q^{(n-m)m}}{[m]_q![n-m]_q!};\,\,[a]_q!=[1]_q
[2]_q\cdots [a]_q;\,\,[a]_q=\frac{q^a-q^{-a}}{q-q^{-1}}$$
From \eqref{10b} one has e.g. $\pp_1(x_2^px_1^k)=(qx_2)^pD_{q^2}x_1^k=(qx_2)^p
[[k]]_{q^2}x_1^{k-1}$ and we consider
\bq\label{14b}
\pp_1(\ga x_1+\gb x_2)^n=\pp_1\sum_0^n\nm_{q^2}(\gb x_2)^m(\ga x_1)^{n-m}=
\end{equation}
$$\sum_0^n\nm_{q^2}(\gb qx_2)^m\ga^{n-m}[[n-m]]_{q^2}x_i^{n-m-1}=$$
$$=\sum_0^{n-1}\nm_{q^2}(\gb qx_2)^m\ga^{n-m}\frac{q^{2(n-m)}-1}{q^2-1}x^{n-m-1}$$
Now from \eqref{13b} one has
\bq\label{15b}
\nm_{q^2}=\frac{[n]_q[n-1]_q!q^{(n-m-1)m}q^m}{[m]_q![n-m]_q[n-m-1]_q!}
\end{equation}
and 
\bq\label{16b}
(\ga x_1+\gb x_2)^{n-1}=\sum_0^{n-1}\left[\begin{array}{c}
n-1\\
m\end{array}\right]_{q^2}(\gb x_2)^m(\ga x_1)^{n-m-1}=\sum^{n-1}_0\Xi_{n,m}
\end{equation}
Consequently
\bq\label{17b}
\pp_1(\ga x_1+\gb x_2)^n=\sum_0^{n-1}\Xi_{n,m}\frac{\ga [n]_qq^{2m}}{[n-m]_q}
\frac{q^{2(n-m)}-1}{q^2-1}=
\end{equation}
$$\frac{\ga
(q^n-q^{-n})}{q^2-1}\sum_0^{n-1}
\Xi_{n,m}q^mq^n\frac{q^{n-m}-q^{m-n}}{q^{n-m}-q^{-n+m}}=
\frac{\ga (q^n-q^{-n})}{q^2-1}\sum_0^{n-1}\Xi_{n,m}q^mq^n$$
$$=\ga[[n]]_{q^2}(\ga x_1+q\gb x_2)^{n-1}$$
Similarly
\bq\label{18b}
\pp_2(\ga x_1+\gb x_2)^n=\pp_2\sum_0^n\nm_{q^2}(\gb x_2)^m(\ga x_1)^{n-m}=
\end{equation}
$$\gb\sum_1^n\nm_{q^2}\frac{q^{2m}-1}{q^2-1}(\gb x_2)^{m-1}(\ga x_1)^{n-m}=$$
$$=\gb\sum _0^{n-1}\frac{[n]_q[n-1]_q!}{[k+1]_q[k]_q!}q^{(n-m)m}\frac{q^{2m}-1}{q^2-1}
(\gb x_2)^k(\ga x_1)^{n-k-1}=$$
$$=\frac{\gb}{q^2-1}\sum_0^{n-1}\Xi_{n,m}\frac{q^n-q^{-n}}{q^{k+1}-q^{-k-1}}
q^{n-k-1}(q^{2(k+1)}-1)=$$
$$=\gb\frac{q^{2n}-1}{q^2-1}(\ga x_1+\gb x_2)^{n-1}=\gb[[n]]_{q^2}
(\ga x_1+\gb x_2)^{n-1}$$
From this one concludes (we will state propositions in a somewhat formal manner)
\begin{proposition}
Let $f(z)=\sum_0^Na_nz^n$ and $D_{q^2}f(z)=f'_{q^2}(z)=\sum_0^N a_nD_{q^2}z^n$.  Then
\bq\label{19b}
f'_{q^2}(z)=\sum_1^Na_n[[n]]_{q^2}z^{n-1}=\sum_0^{N-1}b_mz^m;
\end{equation}
$$\pp_1f(\ga x_1+\gb x_2)=\ga f'_{q^2}(\ga x_1+q\gb x_2);\,\,\pp_2f(\ga x_1+\gb x_2)=
\gb f'_{q^2}(\ga x_1+\gb x_2)$$
\end{proposition}
\indent
From this one computes (setting $g(z)=\sum_0^{N-1}b_mz^m$ in \eqref{19b})
\bq\label{q3}
\pp_1^2f(\ga x_1+\gb x_2)=\gag'_{q^2}(\ga x_1+q\gb x_2);\,\,g'_{q^2}(z)=\sum_1^{N-1}
b_m[[m]]_{q^2}z^{m-1}
\end{equation}
and consequently
\begin{corollary}
For $\ga=\pm\gb$ a solution of the q-wave equation $(\pp_1^2-\pp_2^2)\phi=0$ is given by 
\bq\label{q4}
\phi(x_1,x_2)=F(x_1+x_2)+G(x_1-x_2)
\end{equation}
\end{corollary}
\indent
Next comes questions of uniqueness and boundary conditions of Cauchy type and
surprisingly little seems to have been done in creating a general theory of existence and
uniqueness for q-difference equations, or more generally creating a genuinely algebraic
theory of noncommutative difference equations (see however \cite{fw,kt} for 
some starts and cf. also \cite{ia,lc}).
Before looking at these we want to suggest another approach based on \cite{oh} (cf.
\eqref{4z}-\eqref{8z} and surrounding text).  Working with two variables we have
relations \eqref{1b} and one writes from {\bf (W23)} the relations ${\bf
(W63)}\,\,\mu_1=1+(q^2-1)(x^1\pp_1+x^2\pp_2)$ and $\mu_2=1+(q^2-1)x^2\pp_2$ leading to
{\bf (W24)} and {\bf (W25)}.  Then one gets \eqref{8z} and it seems possible to transport
results about equations in $(X^i,D_i)$ to corresponding results in $(x^i,\pp_i)$.  Recall
$Diff_{GL(n)}\simeq Diff_{q^{-2}}(1)$ and $\simeq$ means isomorphism.  Recall also
\eqref{5z}, or better {\bf (W21)}, to have an isomorphism $Diff(1)\simeq
Diff_{q^{-2}}(1)$.  This seems overly complicated however and one might better go to
another isomorphism between the rings of classical and q-differential operators due to 
Zumino, described in \cite{oh}.  Thus let $x_c^i$ be classical commuting variables and 
choose some ordering of the $x^i$ (noncommutative variables).  Any polynomial $P(x^i)$ is
then written in ordered form and, subsequently, replacing $x^i$ by $x_c^i$ gives
gives a polynomial $\gs(P)(x_c^i)$ determining a ``symbol" map ${\bf (W64)}\,\,\gs:\,{\bf
C}[x^i]\to{\bf C}[x_c^i]$, which is a $\ul{noncanonical}$ isomorphism between the rings
of polynomials.  This provides a map ${\bf (W65)}\,\,\hat{D}\phi=\gs(D(\gs^{-1}(\phi)))$
from the ring of q-difference operators to the ring of classical differential operators
and ${\bf (W66)}\,\,\widehat{D_1D_2}=\hat{D}_1\hat{D}_2$.  The expressions for $\hat{x}^i$
and $\hat{\pp}_i$ will determine the rest.
\\[3mm]\indent
Now suppose one has a solution of a classical wave equation ${\bf
(W67)}\,\,\hat{D}\phi=(\pp_x^2-\pp_y^2)\phi=0$.  This means
$\gs^{-1}(\hat{D}\phi)=D(\gs^{-1}(\phi))$ and we work only on polynomials here with a
fixed ordering for the $x^i$ (while anticipating a possible extension to suitable formal
power series).  Let $\pp_i$ denote classical derivatives with $\pp_i^q$ the corresponding
q-deformed derivatives so 
\bq\label{q5}
\hat{D}\phi=(\pp_x^2-\pp_y^2)\phi(x_c,y_c)\sim\gs\left((\pp_x^q)^2-
(\pp_y^q)^2\right)\hat{\phi}(x,y)
\end{equation}
where $\gs^{-1}\phi\sim\hat{\phi}\sim\phi(x,y)$.  We can use here again \eqref{5z}, or
better {\bf (W21)}, as the definition of $\pp^q_i$ and this leads immediately to
\begin{proposition}
Given a fixed ordering of the $x^i$ and $\gs$ the corresponding isomorphism of {\bf
(W64)}, if $\phi(x_c,y_c)$ is the unique solution of the Cauchy problem $\hat{D}\phi=0$
with polynomial data $\phi(x_c,0)=f(x_c)$ and $(\pp\phi/\pp y_c)(x_c,0)=0$ then
$\hat{\phi}(x,y)=(\gs^{-1}\phi)(x,y)$ is the unique solution of
$D\hat{\phi}=[(\pp_x^q)^2-(\pp_y^q)^2]\hat{\phi}=0$ with $\hat{\phi}(x,0)=f(x)$ and
$\pp_y^q\hat{\phi}(x,0)=0$.
\end{proposition}
\indent
The construction of transmutations from such Cauchy problems is however not immediate
because of possible noncommutativity problems.  For example going back to \eqref{799}
one requires $(P(D_x)-Q(D_y))P(D_x)\phi=P(D_x)(P(D_x)-Q(D_y))\phi=0$ but $PQ\ne QP$ in
general for a noncommutative situation where e.g.
$x_ix_j=qx_jx_i,\,\,\pp_i\pp_j=q^{-1}\pp_j\pp_i$ for $i<j$ as in $\gG_{+}$ (cf.
\eqref{8b} and note that {\bf (W66)} does not stipulate $\widehat{D_1D_2}=\widehat
{D_2D_1}$).  Let us think of $\pp_i=D_{q^2}f$ as in Example 5.1 and use
an ordered $\phi(y,x)=\sum a_{nm}y^nx^m$.  Then \eqref{10b} applies, for $x_1=x$ and 
$x_2=y$ with $xy=qyx,\,\,\pp_y\pp_x=q\pp_x\pp_y,\,\,\pp_xy=qy\pp_x,\,\,\pp_yx=qx\pp_y,\,\,
\pp_xx=q^2x\pp_x+1+(q^2-1)y\pp_y,$ and $\pp_yy=1+q^2y\pp_y$.  One gets
\bq\label{qq5}
\pp_y\phi=\sum a_{nm}(D_{q^2}y^n)x^m=\sum a_{nm}[[n]]_{q^2}y^{n-1}x^m;
\end{equation}
$$\pp_x\phi=\sum a_{nm}(qy)^nD_{q^2}x^m=\sum a_{nm}(qy)^n[[m]]_{q^2}x^{m-1};$$
$$\pp_x\pp_y\phi=\sum a_{nm}[[n]]_{q^2}[[m]]_{q^2}(qy)^{n-1}x^{m-1};\,\,
\pp_y\pp_x\phi=\sum a_{nm}[[n]]_{q^2}[[m]]_{q^2}q^ny^{n-1}x^{m-1}$$
which means ${\bf (W68)}\,\,\pp_y\pp_x\phi=q\pp_x\pp_y\phi$ as could be anticipated from
the 
$\gG_{+}$ rules.  Consequently for $P,Q$ with constant coefficients at least one has
\bq\label{qq6}
P(\pp_x)Q(\pp_y)=(\sum p_n\pp_x^n)(\sum b_m\pp_y^m)=\sum p_nb_m\pp_x^n\pp_y^m;
\end{equation}
$$Q(\pp_y)P(\pp_x)=\sum b_mp_n\pp_y^m\pp_x^n=\sum
b_mp_n\pp_x^n(q\pp_y)^m=P(\pp_x)Q(q\pp_y)$$
Thus, for $P,Q$ with constant coefficients and $\phi=\sum a_{nm}y^nx^m$ ordered, one can
write
\bq\label{qq7}
[P(\pp_x)-Q(\pp_y)]\phi(y,x)=0;\,\,\phi(0,x)=f(x);\,\,\phi(q^{-1}y,0)=Bf(q^{-1}y)
\end{equation}
Set then $\psi(y,x)=P(\pp_x)\phi(q^{-1}y,x)$ so for $y=qz$
\bq\label{qq8}
[P(\pp_x)-Q(\pp_y)]\psi=P(\pp_x)[P(\pp_x)-Q(q\pp_y)]\phi(q^{-1}y,x)=
\end{equation}
$$=P(\pp_x)[P(\pp_x)- Q(\pp_z)]\phi(z,x)=0$$
Hence $\psi(z,0)=B\psi(0,x)=BPf(x)=Q(\pp_z)\psi(z,0)=Q(\pp_z)Bf(z)=Q(q\pp_y)Bf(q^{-1}y)$.
This leads to
\begin{corollary}
It is natural now to (uniquely) define a transmutation operator $B$ (for this fixed
ordering and power series $\phi,\,f$) as $Bf(y)=\phi(0,y)$ (as in \eqref{799})
with $(BPf)(z)=(QBf)(z)$.
\end{corollary}
One can extend the reasoning above to any well posed prolem for a classical PDE (i.e.
well posed means there exists a unique solution) arising from the construction {\bf (W65)}
with suitable boundary conditions (perhaps on vertical and/or horizontal lines) and assert
heuristically
\begin{proposition}
Given a well posed classical problem for $\hat{D}\phi=0$ with polynomial data 
(and solution) and a fixed ordering of noncommutative $x^i$ one obtains a well posed
q-differential problem for $D(\gs^{-1}(\phi))=0$ from {\bf (W65)} and the isomorphism
$\,\,\gs$.
\end{proposition}
\indent
In \cite{fw} one considers wave equations ${\bf (W69)}\,\,[(D_t^{+})^2-{\mc D}_1^{-}{\mc
D}_2^{-}]\phi=0$ where $D_t^{+}\sim t^{-1}(1-T_q)$ and ${\mc D}_x^{-1}\sim x^{-1}
(1-T_q^{-2})$ corresponds to the classical equation $(\pp_t^2-4\pp_1\pp_2)\phi=0$ when
$t\to (1-q)t$ and $x_i\to (1-q^{-1})x_i/2$ with $q\to 1^{-}$.  Solutions to {\bf (W69)}
in q-exponentials have the form ($\ga q^2=\gb\gag$)
\bq\label{q6}
\phi(t,x_1,x_2,\ga,\gb,\gag)=e_q(\ga t)E_q(\gb x_1)E_q(\gag x_2)
\end{equation}
from which one can determine symmetry operators (cf. \cite{fw}).  On the other hand
for ${\bf (W70)}\,\,[D^{-}_t-D_1^{+}D_2^{+}]\phi=0$ one has solutions $\phi=E_q(\ga
t)e_q(\gb x_1)e_q(\gag x_2)$ with $\ga+q\gb\gag=0$.  We recall here from {\bf (W57)} or
\eqref{q2} that $D^{+}_ze_q(\ga z)=\ga e_q(\ga z)$ and $D^{-}_ze_q(\gb z)=-\gb
q^{-1}e_q(\gb z)$ with similar formulas for $E_q$.  For the Helmholz equation
${\bf (W71)}\,\,[D_1^{+}D_2^{+}-\go^2]\phi(x_1,x_2)=0$ solutions can be written in terms
of little q-exponentials via ($\ga\gb=\go^2$)
\bq\label{q7}
\phi(x_1,x_2,\ga,\gb)=e_q(\ga x_1)e_q(\gb x_2)
\end{equation}
For the heat equation in $x,t$ of the form ${\bf (W72)}\,\,[{\mc
D}_t^{-}-(D_x^{+})^2]\phi=0$ one has solutions ${\bf
(W73)}\,\,\phi(t,x,\ga,\gb)=E_{q^2}(\ga t)e_q(\gb x)$ with $\ga+q^2\gb^2=0\,\,(\ga,\gb\in
{\bf C}$).  One arrives at solutions to all these equations be separating variables
according to symmetry operators and their eigenfunctions and this leads for the heat
equation also to solutions
\bq\label{q8}
\phi_n(t,x)=q^{-n(n-3)/2}t^{n/2}H_n\left(\frac{x}{q\sqrt{t}};q\right)
\end{equation}
where $H_n\sim$ discrete q-Hermite polynomial.
\begin{example}
At this point we want to mention some work of Klimek \cite{kt} involving difference
equations and conservation laws (more on this will appear in \cite{ca}).  This
seems to be the only totally partial differential equation 
type material available in a quantum
context and the techniques seem to be natural and powerful.
\end{example}

\subsection{{\bf Kernels and q-special functions.}}

To spell out the nature of the ``transmutation kernels" in \cite{kf} we go to
\eqref{52z}-\eqref{53z} and the preceding table {\bf (A)}-{\bf (E)}.  One sees that
\eqref{52z} could be written as
\bq\label{q35} 
\gb(y,x)=<\gO^P_{\gl}(x),\phi^Q_{\gl}(y)>_{\nu}\sim {\bf Q}(\gO^P_{\gl}(x))\sim
\tl{{\mc P}}(\phi^Q_{\gl}(y));
\end{equation}
$$\gag(x,y)=<\phi^P_{\gl}(x),\gO^Q_{\gl}(y)>_{\go}\sim {\bf P}(\gO^Q_{\gl}(y))=\tl{{\mc
Q}}(\phi_{\gl}^P(x))$$
and it is in these forms that kernels such as $P_{k\ell}$ in \eqref{25z} are displayed.
Some of the other formulas in Section 2 involving Abel and Weyl transforms also stem from
formulas for the calssical situation such as those indicated in Section 3.  We do not try
here to deal with GL or Marchenko ideas in the q-theory but suggest that the framework in
Section 3 seems rich enough to produce analogues in this direction.
\\[3mm]\indent
We make a few further comments here about connecting q-calculus with special functions.
First as an example consider $L=L^{(a,b)}$ and $\mf{L}^{(a,b)}$ of
\eqref{13z}-\eqref{14z} in derivative notation.  Note $T_q^{-1}\sim T_{q^{-1}}$ and write
\bq\label{q32}
\mf{L}=\frac{a}{2}\left(\frac{T_q-1}{x}\right)+\frac{q}{2b}\left(\frac{T_{q^{-1}}-1}{x}
\right)+\frac{1}{2}\left(aT_q+\frac{1}{a}T_{q^{-1}}\right)
\end{equation}
Then recall $T_qf(x)=f(qx)$ with
\bq\label{q33}
D_z^{+}f(z)=\frac{(1-T_q)f(z)}{z};\,\,D_z^{-}f(z)=\frac{(1-T^{-1}_q)f(z)}{z};
\end{equation}
$$D_qf(x)=\frac{(T_q-1)f(x)}{(q-1)x};\,\,\pp_qf(x)=D_{q^2}f(x)$$
(the last equation for the quantum plane).  Then one could write e.g. $T_q\sim
1-zD^{+}_z$ and $T^{-1}_q\sim 1-zD_z^{-}$ with
\bq\label{q34}
\mf{L}=\frac{q}{2}D_x^{+}+\frac{q}{2b}D_x^{-}-\frac{x}{2}\left(aD^{+}_x+
\frac{1}{a}D_x^{-}\right)+(a+a^{-1})
\end{equation}
$$=\frac{a}{2}(1-x)D_x^{+}+\frac{1}{2}
\left(\frac{q}{b}-\frac{x}{a}\right)D_x^{-}+(a+a^{-1})$$
One recalls also $\mf{L}\phi_{\gl}=\gl\phi_{\gl}$ where
$\phi_{\gl}(x;a,b;q)={}_2\phi_1(a\gs,(a/\gs);ab;q;-(bx/z))$ with $\gl=(1/2)(\gs+\gs^{-1})$
(cf. \eqref{12z}).  Concerning a PDE framework in the quantum plane for the eigenfunction
pairings indicated in Section 3 we fix first a suitable integration for
$\pp_x^q,\,\pp_y^q$ operators in \eqref{1b}, \eqref{8b}, or $\gG_{\pm}$ of Example 6.1.
Since from \eqref{39b} one has $\pp_xf(x)\sim D_{q^2}f(x)$ we should require an integral
$\int d_{q^2}x$.  The necessary results are given already in Section 5 but we note here
for background that if $F(y)=\int^y_0f(x)d_{q^2}x$ we want $D_{q^2}F(y)=f(y)$ or
\bq\label{q36}
\frac{F(q^2y)-F(y)}{(q^2-1)y}=f(y)\Rightarrow F(y)=F(q^2y)+(1-q^2)yf(y)\Rightarrow
\end{equation}
$$\Rightarrow F(q^{2n}y)=F(q^{2n+2}y)+(1-q^2)y^{2n}f(q^{2n}y)\Rightarrow
F(y)=F(q^{2n+2}y)+(1-q^2)y\sum_0^nq^{2j}f(q^{2j}y)$$
Then for $|q|<1$ with $F(0)=0$ we have $F(q^{2n+2}y)\to F(0)=0$ and there results
\bq\label{q37}
F(y)=(1-q^2)y\sum_0^{\infty}q^{2j}f(q^{2j}y)
\end{equation}
This gives us a Jackson type integral $\int d_{q^2}x$.
For integrals on infinite intervals one has (cf. \cite{cg,cm,chh,kb,ki})
\bq\label{q38}
\int_0^{\infty(y)}f(x)d_qx=(1-q)\sum_{-\infty}^{\infty}f(q^ky)q^ky;
\end{equation}
$$\int_{-\infty(y)}^{\infty(y)}f(x)d_qx=(1-q)
\sum_{-\infty}^{\infty}[f(q^ky)+f(-q^ky)]q^ky$$
with obvious counterparts for $\int d_{q^2}x$ as in \eqref{q37}.
\\[3mm]\indent
Now we have a nice $q^2$ distribution theory in Section 5 so take here $cos_{q^2}(\mu
x)$ or $Cos_{q^2}(\mu x)$ from \eqref{3e}-\eqref{5e}. From \eqref{10e}
$\pp_q^2z^n\sim (\pp_z^q)^2z^n\sim
D_{q^2}^2z^n=(q^2;q^2)_nz^{n-2}/(q^2;q^2)_{n-2}(1-q^2)^2$ and this is the same as
$[[n]]_{q^2}[[n-1]]_{q^2}z^{n-2}=[(q^{2n}-1)(q^{2n-2}-1)/(q^2-1)^2]z^{n-2}$.
Note from \eqref{2e}-\eqref{3e} or \cite{ka}
\bq\label{q41}
cos_{q^2}(\mu x)=\sum_0^{\infty}\frac{(-1)^k(\mu x)^{2k}}{(q^2;q^2)_{2k}}
\end{equation}
and $[2m]!=q^{-2m^2+m}(q^2;q^2)_{2m}(q^2-1)^{2m}$ indicating a discrepancy with 
\eqref{14a} and \eqref{q39} (as signalled earlier).  In fact \eqref{14a} and \eqref{q39}
correspond to
\bq\label{q411}
Cos_{q^2}(z)=\sum_0^{\infty}\frac{(-1)^nz^{2n}q^{2n(2n-1)}}{(q^2;q^2)_{2n}}
\end{equation}
We will use \eqref{q41} and compute
\bq\label{q42}
(\pp_x^q)^2cos_{q^2}(\mu x)=(\pp_x^q)^2\sum_0^{\infty}\frac{(-1)^k(\mu
x)^{2k}}{(q^2;q^2)_{2k}}=
\end{equation}
$$=\sum_1^{\infty}\frac{(-1)^k\mu^{2k}x^{2k-2}[[2k]]_{q^2}[[2k-1]]_{q^2}}{(q^2;q^2)_{2k}}
=\frac{\mu^2}{(q^2-1)^2}\sum_0^{\infty}\frac{(-1)^m(x\mu)^{2m}}{(q^2;q^2)_{2m}}=
\frac{\mu^2}{(q^2-1)^2}cos_{q^2}(\mu x)$$
since $[[m]]_{q^2}=(q^{2m}-1)(q^2-1)^{-1}=(1-q^{2m})(1-q^2)^{-1}$ and
$(q^2;q^2)_m=\prod_1^m(1-q^{2j})$. 
(We will see below that $\pp_y^q$ needs to be modified for certain situations in the form
e.g. $\mf{D}_{q^2}$ as in {\bf (W83)} below.) 
Now take some second order difference operator
$Q(\pp_y^q)$ (involving $q^2$ differences) with eigenfunctions $\phi^Q_{\mu}(y)$ such that
$Q(\pp_y^q)\phi_{\mu}^Q(y)=-(\mu^2/(q^2-1)^2)\phi^Q_{\mu}(y)$.  Here $\phi^Q_{\mu}(y)$
could involve $q^2$ in various manners.  This is fairly general since if e.g.
$Q(\pp_y^q)\psi_{\nu}(y)=\nu\psi_{\nu}^Q(y)$ then set $\nu=-\mu/(q^2-1)^2$ and rename
$\phi_{\mu}^Q(y)\sim \psi_{-\mu/(q^2-1)^2}^Q(y)$.  We require now further that the
classical corresponding $Q(\pp_y)$ be of say elliptic (symmetric) type so that the Cauchy
problem the classical equation $\pp_x^2\phi=Q(\pp_y)\phi$ will have unique solutions.
Then as indicated in Corollary 6.2 we should be able to define a $q^2$-transmutation
$(\pp_x^q)^2\to Q(\pp_y^q)$ and the kernel should be expressible in the form of a pairing
$<\phi_{\mu}^Q(y),Trig_{\mu}(x)>$ where $Trig_{\mu}$ denotes the appropriate item from
the $q^2$ cosine theory (note the eigenvalue may have to be changed
from \eqref{q42}).  Evidently if we have transmutations $B_1:\,(\pp_x^q)^2\to
Q(\pp_y^q)$ and $B_2: \,(\pp_x^q)^2\to P(\pp_y^q)$ then a transmutation from 
$P(\pp_x^q)\to Q(\pp_y^q)$ can be obtained by the composition $B_1\ci B_2^{-1}$.
\\[3mm]\indent
Thus let us go to Section 5 again and determine the $q^2$ version of the cosine transform
and inversion formulas.  Thus from \eqref{29e}-\eqref{30e} one has
\bq\label{q43}
\mf{F}_{q^2}\phi=\int d_{q^2}z\phi(z){}_0\Phi_1(-;0;q^2;i(1-q^2)q^2zs);
\end{equation}
$$\mf{F}_{q^2}^{-1}\psi=\frac{1}{2\gT_0}\int E_{q^2}(-i(1-q^2)zs)d_{q^2}s\psi(s)$$
where ${\bf (W80)}\,\,{}_0\Phi_1=\ddg E_{q^2}(i(1-q^2)q^qzs)\ddg$.  Then recall
\eqref{q35} where ${\bf (W81)}\,\,\gb(y,x)\sim <\gO_{\gl}^P(x),\phi_{\gl}^Q(y)>_{\nu}\sim
\tl{{\mc P}}(\phi_{\gl}^Q(y))$ and $\tl{{\mc P}}={\mc P}^{-1}\sim
<F(\gl),\gO_{\gl}^P(x)>_{\nu}$.  Since $Cos_{q^2}(z)=(1/2)[E_{q^2}(iz)+E_{q^2}(-iz)]$ one
sees that 
\bq\label{q44}
\mf{F}_{q^2}^{Cos}\phi=\frac{1}{2}\int d_{q^2}z\phi(z)[\ddg E_{q^2}(i(1-q^2)q^2zs)\ddg+
\ddg E_{q^2}(-i(1-q^2)q^2zs)\ddg]=
\end{equation}
$$=\int d_{q^2}z\phi(z)\ddg Cos_{q^2}(1-q^2)q^qzs)\ddg=\int d_{q^2}z\phi(z)
\sum_0^{\infty}\frac{q^{2n(2n-1)}(-1)^n(1-q^2)^{2n}q^{4n}z^{2n}s^{2n}}{(q^2;q^2)_{2n}}$$
The inverse should accordingly be 
\bq\label{q45}
(\mf{F}_{q^2}^{Cos})^{-1}=\frac{1}{2\gT_0}\int d_{q^2}s\psi(s)\ddg Cos(1-q^2)zs)\ddg=
\end{equation}
$$=\frac{1}{2\gT_0}\int d_{q^2}s\psi(s)\sum_0^{\infty}\frac{(-1)^nq^{2n(2n-1)}(1-q^2)^{2n}
z^{2n}s^{2n}}{(q^2;q^2)_{2n}}$$
The fact that this is true follows from \eqref{37e}.
Note that one could also develop a Fourier cosine theory following
\eqref{14a}-\eqref{16a} but the approach in section 5 from \cite{oe} seems much better;
in particular it is more complete with useful formulas like \eqref{38e} for example as
well as the distribution format.  In any event we could possibly think here of ${\mc P}$
as the $Cos_{q^2}$ transform $\mf{F}_{q^2}$ with $\phi_{\gl}^P\sim \ddg Cos_{q^2}((1-q^2)
q^2zs)\ddg$.  This has to be slightly modified however since (cf. \eqref{q42}) ${\bf
(W82)}
\,\,\pp_z^2z^{2n}=[[2n]]_{q^2}[[2n-1]]_{q^2}z^{2n-2}=(q^{4n}-1)(q^{4n-4}-1)/(q^2-1)^2
z^{2n-2}$
with
\bq\label{q47}
\pp_z^2\ddg
Cos_{q^2}((1-q^2)q^2zs)\ddg=
\end{equation}
$$=\sum_1^{\infty}\frac{q^{2n(2n-1)}(-1)^n(1-q^2)^{2n}
q^{4n}z^{2n-2}s^{2n}(q^{4n}-1)(q^{4n-4}-1)}{(q^2;q^2)_{2n}(q^2-1)^2}=$$
$$=-s^2q^6\sum_0^{\infty}\frac{(-1)^mq^{2m(2m-1)}q^{4m}(q^4z)^{2m}s^{2m}}{(q^2;q^2)_
{2m}}=-s^2q^6\ddg Cos_{q^2}(1-q^2)q^2(q^4z)s\ddg$$
This shows that $(\pp_z^q)^2$ 
is not the right operator to use when looking for
eigenfunctions.  We would rather have an operator such that the argument of
$Cos_{q^2}(1-q^2)q^2zs)$ can be maintained.  Thus we need a power of $q$ to offset the
$q^{2n(2n-1)}$ term or say $q^{-8m}=q^{-8(n-1)}=q^{-8n+8}$ so define ${\bf (W83)}\,\,
\mf{D}_{q^2}z^n=[[n]]_{q^2}z^{n-1}(q^2)^{\ga_n}$ with $\ga_{2n}+\ga_{2n-1}=-4n+\gb$ for
$\gb$ any constant.  Then take e.g. $\ga_n=-n+\gag$ so
$\ga_{2n}+\ga_{2n-1}=-2n-2n+1+2\gag=-4n+(2\gag+1)$ and one can remove the $q^6$ term as
well if $2(-4n+2\gag+1)=-8n+8-6=-8n+2$.  Taking $\gag=0$ one arrives at
${\bf (W84)}\,\,\mf{D}_{q^2}z^n=[[n]]_{q^2}z^{n-1}(q^2)^{-n}$ and hence
\bq\label{q48}
\mf{D}_{q^2}^2\ddg Cos_{q^2}((1-q^2)q^2zs)\ddg=
\end{equation}
$$=-s^2\sum_0^{\infty}\frac{(-1)^m
(1-q^2)^{2m}q^{4m}z^{2m}s^{2m}}{(q^2;q^2)_{2m}}=-s^2\ddg Cos_{q^2}((1-q^2)q^2zs)\ddg$$
\begin{proposition}
Take $\phi_s^P(z)\sim\ddg Cos_{q^2}((1-q^2)q^2zs)\ddg$ with $P\sim \mf{D}^2_{q^2}$ and
$P_z\phi^P_s(z)=-s^2\phi_s^P(z)$ while ${\mc P}\phi(s)=(\mf{F}_{q^2}^{Cos}\phi)(s)$.  The
kernel in {\bf (W81)} will then take the form
\bq\label{q49}
\gb(y,z)=\tl{{\mc P}}(\phi_{\gl}^Q(y))=(\mf{F}_{q^2}^{Cos})^{-1}(\phi_s^Q(y))=
\frac{1}{2\gT_0}\int d_{q^2}s\phi_s^Q(y)\ddg Cos((1-q^2)zs)\ddg
\end{equation}
which is represented by a series as in Section 5.
\end{proposition}
\indent
This provides a beginning for implementation of the program indicated in Section 1.
There are a number of obvious matters to investigate now besides looking at examples and
convergence questions.  For example connections to integrable systems and tau functions
have quantum counterparts and the GLM machinery is related to this.  Relations between
classical generating functions and tau functions should be further explored in the spirit
mentioned in Section 1.  Further devleopment of quantum calculi with applications to
differential and diffference equations seems inevitable. Generally the fundamental
combinatorial material arising via braiding \`a la R matrices or quantum groups has in
fact defined the nature of quantum algebra-calculi and probably some diversity in this
and related physical science will persist.  We will return to this in \cite{ca}.

\end{document}